\documentclass[12pt]{amsart}

\usepackage{amsmath}
\usepackage{amssymb}
\usepackage{caption}
\usepackage{subcaption}
\usepackage{float}
\usepackage[curve]{xypic}
\usepackage{cite}
\usepackage{enumitem}
\usepackage{color}
\usepackage{url}
\usepackage[usenames,dvipsnames]{xcolor}
\usepackage{graphicx}
\usepackage{verbatim}
\usepackage{tikz}
\usetikzlibrary{arrows}
\usepackage[normalem]{ulem}
\usepackage{hyperref}
\hypersetup{
     colorlinks   = true,
     citecolor    = black,
     linkcolor    = blue
}

\makeatletter 
\def\@cite#1#2{{\m@th\upshape\bfseries%
[{#1\if@tempswa{\m@th\upshape\mdseries, #2}\fi}]}}
\makeatother 

\theoremstyle{plain}
\newtheorem{thm}{Theorem}[section]

\newtheorem{cor}[thm]{Corollary}

\theoremstyle{definition}

\newtheorem{ex}[thm]{Example}

\newtheorem{prob}[thm]{Problem}
\newtheorem{conj}[thm]{Conjecture}
\theoremstyle{}

\numberwithin{equation}{subsection}
\renewcommand{\bold}[1]{\medskip \noindent {\bf #1 }\nopagebreak}

\newcommand{\nc}{\newcommand}
\newcommand{\rnc}{\renewcommand}
\newcommand{\e}{\varepsilon}

\nc\bA{\mathbb{A}}
\nc\bB{\mathbb{B}}
\nc\bC{\mathbb{C}}
\nc\bD{\mathbb{D}}
\nc\bE{\mathbb{E}}
\nc\bF{\mathbb{F}}
\nc\bG{\mathbb{G}}
\nc\bH{\mathbb{H}}
\nc\bI{\mathbb{I}}
\nc{\bJ}{\mathbb{J}} 
\nc\bK{\mathbb{K}}
\nc\bL{\mathbb{L}}
\nc\bM{\mathbb{M}}
\nc\bN{\mathbb{N}}
\nc\bO{\mathbb{O}}
\nc\bP{\mathbb{P}}
\nc\bQ{\mathbb{Q}}
\nc\bR{\mathbb{R}}
\nc\bS{\mathbb{S}}
\nc\bT{\mathbb{T}}
\nc\bU{\mathbb{U}}
\nc\bV{\mathbb{V}}
\nc\bW{\mathbb{W}}
\nc\bY{\mathbb{Y}}
\nc\bX{\mathbb{X}}
\nc\bZ{\mathbb{Z}}
\nc\cA{\mathcal{A}}
\nc\cB{\mathcal{B}}
\nc\cC{\mathcal{C}}
\rnc\cD{\mathcal{D}}
\nc\cE{\mathcal{E}}
\nc\cF{\mathcal{F}}
\nc\cG{\mathcal{G}}
\rnc\cH{\mathcal{H}}
\nc\cI{\mathcal{I}}
\nc{\cJ}{\mathcal{J}} 
\nc\cK{\mathcal{K}}
\rnc\cL{\mathcal{L}}
\nc\cM{\mathcal{M}}
\nc\cN{\mathcal{N}}
\nc\cO{\mathcal{O}}
\nc\cP{\mathcal{P}}
\nc\cQ{\mathcal{Q}}
\rnc\cR{\mathcal{R}}
\nc\cS{\mathcal{S}}
\nc\cT{\mathcal{T}}
\nc\cU{\mathcal{U}}
\nc\cV{\mathcal{V}}
\nc\cW{\mathcal{W}}
\nc\cY{\mathcal{Y}}
\nc\cX{\mathcal{X}}
\nc\cZ{\mathcal{Z}}
\nc\wt{\widetilde}

\newcommand{\bk}{{\mathbf{k}}}

\newcommand{\SL}{\mathrm{SL}}
\newcommand{\GL}{\mathrm{GL}}

\renewcommand{\GL}{GL^+(2,\bR)}

\newcommand{\wrt}[1]{\mathrm{d}{#1}}

\newcommand{\IC}{\iota}

\nc{\dmo}{\DeclareMathOperator}
\rnc{\Re}{\operatorname{Re}}
\rnc{\Im}{\operatorname{Im}}
\dmo{\rank}{rank}
\dmo{\End}{End}
\dmo{\Hom}{Hom}
\dmo{\Jac}{Jac}
\dmo{\Id}{Id}
\dmo{\Ann}{Ann}
\dmo{\Area}{Area}
\dmo{\CP}{\bC P^1}
\dmo{\Aut}{Aut}
\dmo{\MCG}{MCG}
\dmo{\Vol}{Vol}
\dmo{\dVol}{\wrt{Vol}}
\dmo{\dd}{d}
\dmo{\Homeo}{Homeo}
\dmo{\Stab}{Stab}
\dmo{\ML}{\mathcal{ML}}
\dmo{\MF}{\mathcal{MF}}
\dmo{\Prob}{Prob}
\dmo{\QD}{QD}
\dmo{\Tw}{Tw}

\title[Mirzakhani's work on Riemann surfaces]{A tour through Mirzakhani's work on moduli spaces of Riemann surfaces}
\author[Wright]{Alex~Wright}
\AtBeginDocument{%
   \def\MR#1{}
}

\begin{document}
\vspace{-.5cm}
\maketitle
\thispagestyle{empty}


\setcounter{tocdepth}{1} 
\vspace{-.7cm}
\renewcommand\contentsname{} 

\begingroup
\let\clearpage\relax
\vspace{-1cm} 
\tableofcontents
\endgroup
\vspace{-1.2cm}
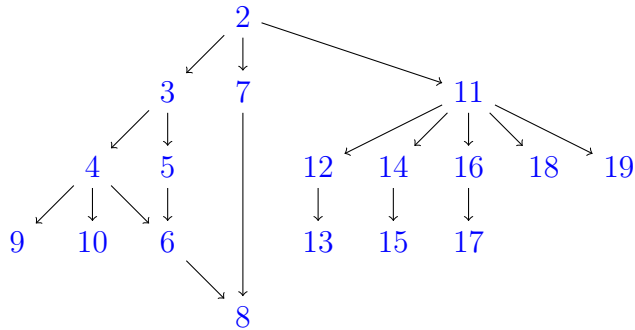
\begin{figure}[ht!]
\begin{tikzpicture}
    \node (PreTeich) at (0,0) {\ref{S:PreTeich}};
    
    \node (VolSymp) at (0,-1) {\ref{S:VolSymp}};
    \node (VolM11) at (-1,-1) {\ref{S:VolM11}};
    \node (PreTD) at (3,-1) {\ref{S:PreTD}};
   \draw [->] (PreTeich) -- (VolM11);
   \draw [->] (PreTeich) -- (VolSymp);
   \draw [->] (PreTeich) -- (PreTD);
   
   \node (McShane) at (-1,-2) {\ref{S:McShane}};
   \node (Integrating) at (-2,-2) {\ref{S:Integrating}};
   \draw [->] (VolM11) -- (McShane);
   \draw [->] (VolM11) -- (Integrating);
   
   \node (Recursive) at (-1,-3) {\ref{S:Recursive}};
   \draw [->]  (McShane) -- (Recursive);
   \draw [->]  (Integrating) -- (Recursive);
   
   \node (Witten) at (-0,-4) {\ref{S:Witten}};
   \draw [->]  (Recursive) -- (Witten);
   \draw [->]  (VolSymp) -- (Witten);

  \node (CountSimple) at (-3,-3) {\ref{S:CountSimple}};
  \node (Random) at (-2,-3) {\ref{S:Random}};

  \draw [->]  (Integrating) -- (CountSimple);
  \draw [->]  (Integrating) -- (Random);
  
  \node (Earthquake) at (1,-2) {\ref{S:Earthquake}};
  \node (ABEM) at (2,-2) {\ref{S:ABEM}};
  \node (Measures) at (3,-2) {\ref{S:Measures}};
  
  \draw [->]  (PreTD) -- (Earthquake);
  \draw [->]  (PreTD) -- (ABEM);
  \draw [->]  (PreTD) -- (Measures);
  
  \node (Ergodic) at (1,-3) {\ref{S:Ergodic}};
    \node (Orbits) at (2,-3) {\ref{S:Orbits}};
  \node (OrbitClosures) at (3,-3) {\ref{S:OrbitClosures}};
  
  \draw [->]  (Earthquake) -- (Ergodic);
  \draw [->]  (ABEM) -- (Orbits);
  \draw [->]  (Measures) -- (OrbitClosures);
  
  \node (Effective) at (4,-2) {\ref{S:Effective}};
  \draw [->]  (PreTD) -- (Effective);
  
  \node (RW) at (5,-2) {\ref{S:RW}}; 
  \draw [->]  (PreTD) -- (RW);
  
\end{tikzpicture}
\caption{Section pre-requisite chart. If interested in random structures, aim for \ref{S:Random}; if  the Prime Number Theorem, \ref{S:CountSimple} or \ref{S:Orbits} or \ref{S:Effective}; if algebraic or symplectic geometry, \ref{S:Witten}; if dynamics, anything following \ref{S:PreTD}.}
\label{F:Pre}
\end{figure}
\vfill


\section{Introduction}\label{S:intro}

This survey aims to be a tour through Maryam Mirzakhani's  remarkable work on Riemann surfaces, dynamics, and geometry. The star characters, all across mathematics and physics as well as in this survey, are the moduli spaces $\cM_{g,n}$ of Riemann surfaces. 

Sections \ref{S:PreTeich} through \ref{S:Random} all relate to Mirzakhani's study of the size of these moduli spaces, as measured by the Weil-Petersson symplectic form.  Goldman has shown that  many related moduli spaces also have a Weil-Petersson symplectic form, so this can be viewed as part of a broader story \cite{Goldman:Symplectic}. Even more important than the broader story, Mirzakhani's study unlocks applications to the topology of $\cM_{g,n}$, random surfaces of large genus, and even geodesics on individual hyperbolic surfaces. 

Sections \ref{S:PreTD} to \ref{S:RW} reflect the philosophy that $\cM_{g,n}$, despite being a totally inhomogeneous object, enjoys many of the dynamical  properties of nicer spaces, and even some of the  dynamical miracles characteristic of homogeneous spaces. The dynamics of group actions in turn clarify the geometry of $\cM_{g,n}$ and produce otherwise unattainable counting results.  

Our goal is  not to provide a comprehensive reference, but rather to highlight some of the most beautiful and easily understood ideas from the roughly 20 papers that constitute Mirzakhani's work in this area. Very roughly speaking, we devote comparable time to each paper or closely related group of papers. This means in particular that we cannot proportionately discuss  the longest paper \cite{EskinMirzakhani:InvariantMeasures}, but on this topic the reader may see the surveys \cite{Zorich:Magic, Wright:Billiards, Quint:Rigidite}.  We include some open problems, and hope that we have succeeded in conveying the thriving legacy of Mirzakhani's research. 

We invite the reader to discover for themselves Mirzakhani's five papers on combinatorics \cite{MahmoodianMirzakhani:Tripartite, Mirzakhani:Choosable, Mirzakhani:Shur, Mirzakhani-Vondrak:SpernerFair, MirzakhaniVondrak:SpernerOptimal}, where the author is not qualified to guide the tour. 

We also omit comprehensive citations to work preceding Mirzakhani, suggesting instead that the reader may get off the tour bus at any time to find more details and context in the references and re-board later, or to revisit the tour locations  at a later date. Where possible, we give references to expository sources, which will be more useful to the learner than the originals.  The reader who consults the references  will be rewarded with views of the vast tapestry of important and beautiful work that Mirzakhani builds upon, something we can only offer tiny glimpses of here. 

We hope that a second year graduate student who has previously encountered the definitions of hyperbolic space, Riemann surface, line bundle, symplectic manifold, etc., will be able to read and appreciate the survey, choosing not to be distracted by the occasional remark aimed at the experts. 

Other surveys on Mirzakhani's work include \cite{Wolpert:Families, Wolpert:Mirzakhani, Do:BigSurvey, McMullen:ICM, Zorich:Magique, Zorich:Magic, Huang:Mirzakhani, Wright:Billiards, Quint:Rigidite, Martin:Early,Wright:Earthquake}. See also the issue of the Notices of the AMS that was devoted to Mirzakhani \cite{NoticesMM}.

\bold{Acknowledgments.} We are especially grateful to Scott Wolpert and Peter Zograf for helpfully answering an especially large number of questions, to Paul Apisa for especially detailed comments that resulted in major  improvements, and to Francisco Arana Herrera, Jayadev Atheya, Alex Eskin, Mike Lipnowski, Curt McMullen, Amir Mohammadi, Juan Souto, Kasra Rafi, and Anton Zorich for in depth conversations. 
 
 We also thank the many people who have commented on earlier drafts, including, in addition to the above,
 Giovanni Forni,
Dmitri Gekhtman, %
Mark Greenfield, 
Chen Lei,
Ian Frankel, 
Athanase Papadopoulos,
Bram Petri, 
Feng Zhu, 
and the anonymous referee. 

\section{Preliminaries on Teichm\"uller Theory}\label{S:PreTeich}

We begin with the beautiful and basic results that underlie most of Mirzakhani's work. 

\subsection{Hyperbolic geometry and complex analysis.}  All surfaces are assumed to be orientable and connected.  Any simply connected surface with a complete Riemannian metric of constant curvature $-1$ is isometric to the upper half plane $$\bH=\{x+iy \mid y>0\}$$ endowed with the hyperbolic metric 
$$(ds)^2 = \frac{(dx)^2+(dy)^2}{y^2}.$$
Perhaps the most important miracle of low dimensional geometry is that the group of orientation preserving isometries of hyperbolic space is equal to the group of biholomorphisms of $\bH$. (Both are equal to the group $PSL(2,\bR)$ of M\"obius transformations that stabilize the upper half plane.) 

Every oriented complete hyperbolic surface $X$ has universal cover $\bH$, and the deck group  acts on $\bH$ via orientation preserving isometries. Since these isometries are also biholomorphisms, this endows $X$ with the structure of a Riemann surface, namely an atlas of charts to $\bC$ whose transition functions are biholomorphisms. 

Conversely, every  Riemann surface $X$ that is not simply connected, not $\bC\setminus \{0\}$, and not a torus has universal cover $\bH$, and the deck group acts on $\bH$ via biholomorphisms. Since these biholomorphisms are also  isometries, this endows $X$ with a complete hyperbolic metric. 

\subsection{Cusps, geodesics, and collars.} Suppose $X$ is a complete hyperbolic surface. Each subset of $X$ isometric to 
$$\{z=x+iy\in \bH : y>1\}/\langle z\mapsto z+1\rangle$$
is called a cusp. Each cusp has infinite diameter and finite volume. Distinct cusps are disjoint, and if $X$ has finite area then the complement of the cusps is compact. 

Each cusp is biholomorphic to a punctured disc via the exponential map. If $X$ has finite area, then $X$ is biholomorphic to a compact Riemann surface minus a finite set of punctures, and the punctures are in bijection with the cusps.

Any closed curve on $X$ not homotopic to a point or a loop around a cusp is isotopic to a unique closed geodesic. Unless otherwise stated, all closed curves we consider will be of this type. A closed geodesic is called simple if it does not intersect itself.

\begin{figure}[ht!]
\includegraphics[width=.5\linewidth]{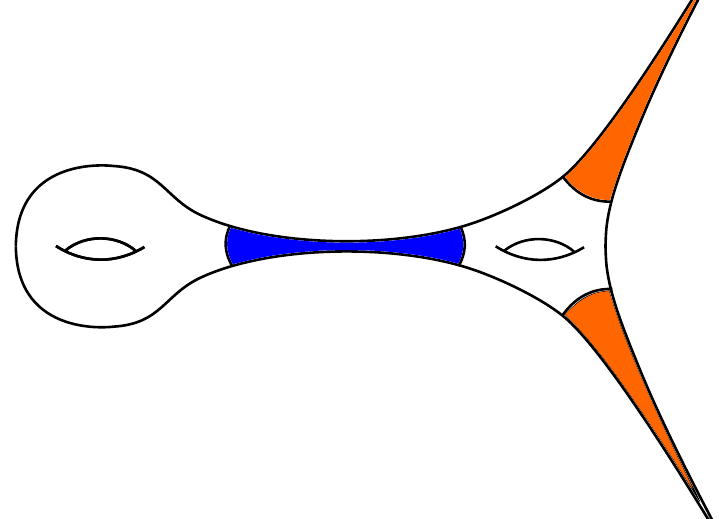}
\caption{Two cusps and a collar.}
\label{F:collars}
\end{figure}

Gauss-Bonnet gives that any closed hyperbolic surface of genus $g$ has area $2\pi(2g-2)$. This is the first indication that a hyperbolic surface cannot be ``small". Moreover, the Collar Lemma gives that any closed geodesic of length less than a universal constant is simple, and every short simple closed geodesic must be surrounded by a large embedded annulus known as its collar. As the length of the simple closed geodesic goes to zero, the size of its collar goes to infinity.

\subsection{Building a surface out of pants.} A significant amount of this survey will concern hyperbolic surfaces with boundary. We will always assume that any surface with boundary that we consider can be isometrically embedded in a complete surface so that the boundary consists of a finite union of closed geodesics. 

A hyperbolic sphere with three  boundary components is known as a pair of pants, or simply as a pants. A fundamental fact gives that, for any three numbers $L_1, L_2, L_3>0$, there is a unique  pants with these three boundary lengths.  Each $L_i$ may also be allowed to be zero, in which case a pants, now somewhat degenerate, has a cusp instead of a boundary component.  

\begin{figure}[ht!]
\includegraphics[width=.35\linewidth]{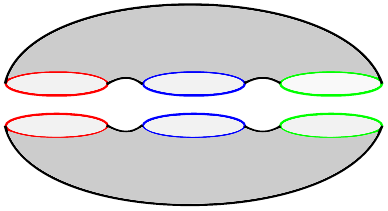}
\caption{Gluing pants to form a genus two surface.}
\label{F:TwoPants}
\end{figure}

One of the simplest ways to build a closed hyperbolic surface is by gluing together pants. For example, given two pants with the same boundary lengths, we may glue together the corresponding boundaries to obtain a closed genus 2 surface. In fact, the corresponding boundaries can be glued using different isometries from the circle to the circle, giving infinitely  many genus two hyperbolic surfaces. More complicated surfaces can be obtained by gluing together more pants.

\subsection{Teichm\"uller space and moduli space.} We define moduli space $\cM_{g,n}$ formally as the set of equivalence classes of oriented genus $g$ hyperbolic surfaces with $n$ cusps labeled by $\{1, \ldots, n\}$, where two surfaces are considered equivalent if they are isometric via an orientation preserving isometry that respects the labels of the cusps. Equivalently, $\cM_{g,n}$ can be defined as the set of equivalence classes of genus $g$ Riemann surfaces with $n$ punctures labeled by $\{1, \ldots, n\}$, where two surfaces are considered equivalent if they are biholomorphic via a biholomorphism that respects the labels of the punctures.

We will follow the almost universal abuse of referring to a point in $\cM_{g,n}$ as a hyperbolic or Riemann surface, leaving out the notational bookkeeping of the equivalence class. 

Teichm\"uller space $\cT_{g,n}$ is defined to be the set 
$$\cT_{g,n}=\{(X,[\phi])\}$$
of points $X$ in  $\cM_{g,n}$, which as indicated we think of as hyperbolic or Riemann surfaces, equipped with a homotopy class $[\phi]$ of orientation preserving homeomorphisms $\phi:S_{g,n} \to X$ from a fixed oriented topological surface $S_{g,n}$ of genus $g$ with $n$ punctures. The homotopy class $[\phi]$ is called a marking, and one says that $\cT_{g,n}$ parametrizes marked hyperbolic or Riemann surfaces. 

Let $\Homeo^+(S_{g,n})$ denote the group of orientation preserving homeomorphisms of $S_{g,n}$ that do not permute the punctures. This group acts on $\cT_{g,n}$ by precomposition with the marking. The subgroup $\Homeo_0^+(S_{g,n})$ of homeomorphisms isotopic to the identity acts trivially, so the quotient 
$$\MCG_{g,n}= \Homeo^+(S_{g,n}) / \Homeo_0^+(S_{g,n})$$
acts on $\cT_{g,n}$. This countable group is called the mapping class group, and 
$$ \cM_{g,n} = \cT_{g,n}/ \MCG_{g,n}.$$

Given $L=(L_1, \ldots, L_n)\in \bR_+^n$, we can similarly define $\cT_{g,n}(L)$ to be the Teichm\"uller space of oriented genus $g$ hyperbolic surfaces with $n$ boundary components of length $L_1, \ldots, L_n$, and $\cM_{g,n}(L)$ to be the corresponding moduli space. Here $S_{g,n}$ is replaced with a genus $g$ surface with $n$ boundary circles, and we define $\Homeo^+(S_{g,n})$ to be the orientation preserving homeomorphisms that do not permute the boundary components.  

It is sometimes convenient to allow $L_i=0$ in the definitions above, in which case the corresponding boundary is replaced by a cusp. For example, using this convention $\cM_{g,n} = \cM_{g,n}(0, \ldots, 0)$. 

\subsection{Classification of simple closed curves.} Let $\alpha$ and $\beta$ be two different simple closed curves on $S_{g}$ that are non-separating, in that cutting either curve does not disconnect the surface. In this case the result of cutting either $\alpha$ or $\beta$ is homeomorphic to a genus $g-1$ surface with two boundary curves, and hence are homeomorphic to each other. This homeomorphism can be modified to give rise to a homeomorphism of $S_g$ that takes $\alpha$ to $\beta$. In particular, we conclude that there is some $f\in \MCG_g$ such that $f([\alpha])=[\beta]$, where $[\alpha]$ denotes the homotopy  class of $\alpha$. 

Next suppose that $\alpha$ is a separating simple closed curve. In this case, $S_g\setminus \alpha$ has two components, one of which is a surface of genus $g_1$ with one boundary component, and the other of which is a surface of genus $g_2=g-g_1$ with one boundary component. If $\beta$ is another separating curve, then there is some $f\in \MCG_g$ such that $f([\alpha])=[\beta]$ if and only if the set $\{g_1, g_2\}$ arising from $\beta$ is the same as for $\alpha$. 

In summary, there is a single mapping  class group orbit of non-separating simple closed curves on $S_g$, and $\lfloor \frac{g}2 \rfloor$ mapping class group orbits of separating simple closed curves. 

\subsection{The twist flow.}\label{SS:twist} Let $\alpha$ be a simple closed curve on $S_{g,n}$ that isn't a loop around a cusp.  We now introduce the twist flow $\Tw_t^\alpha$ on Teichm\"uller space as follows. It may be conceptually helpful to start by assuming $t$ is small and positive. 

For each point $(X,[\phi])\in \cT_{g,n}$, we can consider the geodesic representative of $\phi(\alpha)$. Cut this geodesic to obtain a surface with two geodesic boundary components of equal length. Both of these components inherit an orientation from the surface, and we will call the positive direction ``left". 
\begin{figure}[ht!]
\includegraphics[width=.35\linewidth]{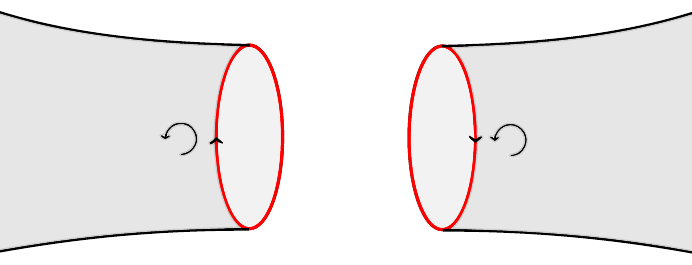}
\caption{The orientation  the  boundary components.}
\label{F:TwistLeft}
\end{figure}
Re-glue the two components by the original identification, composed with a rotation by $t$, so that if two points $p_1, p_2$ on the two boundary components were originally identified, now $p_1$ is identified with the point $t$ to the left of $p_2$, and vice versa.   
\begin{figure}[ht!]
\includegraphics[width=\linewidth]{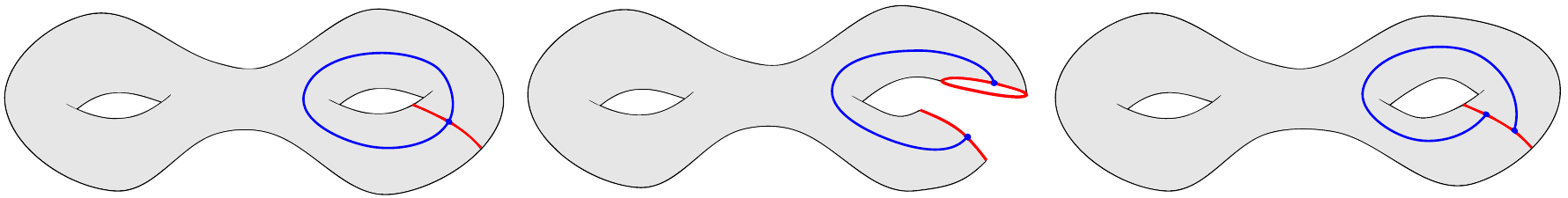}
\caption{The effect of a twist on a transverse curve.}
\label{F:Twist}
\end{figure}

If we use the notation $\Tw_t^\alpha(X,[\phi])=(X_t,[\phi_t])$, this re-gluing defines $X_t$. The marking $[\phi_t]$ is more subtle, and we will omit its definition. Here it will suffice to accept that, despite the fact that $X_{t+\ell_\alpha(X)}=X_t$, the twist path is injective, so $(X_{t_1},[\phi_{t_1}])=(X_{t_2},[\phi_{t_2}])$ if and only if $t_1=t_2$. 

\subsection{Fenchel-Nielsen coordinates.} Fix a pants decomposition of $S_{g,n}$. This is a  collection of disjoint  simple closed curves, such that cutting these curves gives a collection of topological pants. It turns out any such collection has $3g-3+n$ curves, and we denote these curves $\{\alpha_i\}_{i=1}^{3g-3+n}$.  If $n>0$,  some of the pants will be degenerate, in that they will have a puncture instead of boundary circle. 

\begin{figure}[ht!]
\includegraphics[width=.5\linewidth]{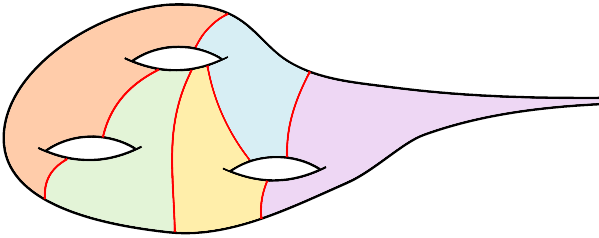}
\caption{A pants decomposition.}
\label{F:Pants}
\end{figure}

Given a marked hyperbolic surface $(X,\phi)$, we can consider the curves $\phi(\alpha_i)$ on $X$. Let $\ell_{\alpha_i}(X)$ denote the length of the geodesic homotopic to $\phi(\alpha_i)$. For short, we write $\ell_i$ to denote $\ell_{\alpha_i}(X)$. 

Each hyperbolic surface in $\cT_{g,n}$ can be obtained by gluing together the pants with the correct boundary lengths in the correct combinatorial pattern, but additional  parameters are required in the construction to keep track of how the boundary curves are glued together.

\begin{thm}[Fenchel-Nielsen]
There are functions $$\tau_i: \cT_{g,n}\to \bR,\quad i=1, \ldots, 3g-3+n$$
such that the map 
$\cT_{g,n} \to \bR_+^{3g-3+n} \times \bR^{3g-3+n}$ defined by $$(\ell_1, \ldots, \ell_{3g-3+n}, \tau_1, \ldots, \tau_{3g-3+n})$$
is a homeomorphism, and so that for each $i$ and all $t$, $$\tau_i(\Tw_{\alpha_i}^t(X,[\phi]))=t+\tau_i(X,[\phi]),$$ and all the other coordinates of $\Tw_{\alpha_i}^t(X,[\phi])$ and $(X,[\phi])$ are the same. 
\end{thm}

The twist parameters $\tau_i$ and the length parameters $\ell_i$ are called Fenchel-Nielsen coordinates for Teichm\"uller space.

One can show that the mapping class group acts properly discontinuously on $\cT_{g,n}$. In particular, the stabilizer of each point is finite. The quotient $\cM_{g,n}$ is thus an orbifold, which is similar to a manifold except that some points have neighborhoods homeomorphic to a neighbourhood of the origin in $\bR^{6g-6+2n}$ quotiented by a finite group action. 

Fenchel-Nielsen coordinates work similarly for $\cT_{g,n}(L)$. Note there are no twist or length parameters for the geodesic boundary curves, since they have fixed lengths $L_i$ and are not glued to anything. 


\subsection{The Weil-Petersson symplectic structure.} Fix a choice of Fenchel-Nielsen coordinates, and define 
$$\omega_{WP} = \sum d\ell_i \wedge d\tau_i$$
to be the standard symplectic form in these coordinates on $\cT_{g,n}(L)$. 

Wolpert proved that this symplectic form is invariant under the action of the mapping class group. Hence it descends to a symplectic form $\omega_{WP}$ on $\cM_{g,n}$ or $\cM_{g,n}(L)$.\footnote{We define $\omega_{WP}$ on moduli space so that its pullback to Teichm\"uller space is the standard symplectic form defined above. In other words, it is the standard symplectic form in local Fenchel-Nielsen coordinates on moduli space.

This is sometimes called the ``topologist's definition", and it ignores that $\cM_{g,n}$ may be considered as a stack, which is the algebro-geometric version of an orbifold.    To reconcile with the algebro-geometric perspective without using stacks, one could also define the Weil-Petersson volume form on $\cM_{g,n}$  as the local push-forward of the Weil-Petersson volume form on Teichm\"uller space. The  definitions give volume forms that are equal except for $\cM_2$ and $\cM_1(L)$, where the reconciled volume form is half of the topologist's. Every surface in those two moduli spaces has an involution symmetry.

Independently of this issue, it is also common to include a separate factor of $\frac12$ in the definition of $\omega_{WP}$ for all $g$ and $n$ \cite[Section 5]{Wolpert:Cusps}.} 

Wolpert also showed that in the case of $\cT_{g,n}$, this symplectic form is twice the one arising from the   Weil-Petersson K\"ahler structure on $\cT_{g,n}$. This result is sometimes called Wolpert's Magic Formula, since 
\begin{itemize} 
\item the definition of the Weil-Petersson K\"ahler structure, although very natural, gives no hint of a relationship to Fenchel-Nielsen coordinates, and 
\item it is surprising that the two-form $\sum d\ell_i \wedge d\tau_i$, obtained from a pants decomposition, does not depend on which pants decomposition is used. 
\end{itemize} 

The associated Weil-Petersson volume form, which is the standard volume form in local Fenchel-Nielsen coordinates, is the most natural known notion of volume on moduli space. The Weil-Peterson volume of each moduli space is finite. 


\subsection{References.} More details can be found in the books \cite[Chapters 10, 12]{FarbMargalit:Primer} and, for the Weil-Petersson symplectic structure, \cite[Chapter 3]{Wolpert:Families}.

\section{The volume of $\cM_{1,1}$}\label{S:VolM11}

Mirzakhani discovered an elegant new computation of the volume of $\cM_{1,1}$. We reproduce this computation, which is highlighted in the introduction to \cite{Mirzakhani:Invent}, since it was perhaps the first seed for her thesis. The starting point is the  remarkable identity 
$$\sum_\alpha \frac{1}{1+e^{\ell_\alpha(X)}}= \frac12,$$
of McShane \cite{McShane:Identity}, which gives that a certain sum involving the lengths $\ell_\alpha(X)$ of all simple closed geodesics $\alpha$ on $X\in \cM_{1,1}$ is independent of  $X\in \cM_{1,1}$. In Section \ref{S:McShane} we'll explain where this identity comes from. 

Let $\cM_{1,1}^*$ denote the infinite cover of $\cM_{1,1}$ parametrizing pairs $(X, \alpha)$, where $X\in \cM_{1,1}$ and $\alpha$ is a simple closed geodesic on $X$. Mirzakhani's computation begins 
 
\begin{eqnarray*}
\frac12 \Vol(\cM_{1,1}) &=& \int_{\cM_{1,1}} \sum_\alpha \frac{1}{1+e^{\ell_X(\alpha)}} \dVol_{WP}
\\&=& \int_{\cM_{1,1}^*}  \frac{1}{1+e^{\ell_X(\alpha)}} \dVol_{WP}. 
\end{eqnarray*}
This ``unfolding" is justified because the fibers of the map $\cM_{1,1}^*\to \cM_1$ are precisely the set of simple closed geodesics, which is the set being summed over in McShane's identity. 

Given a point $(X,\alpha)\in \cM_{1,1}^*$, we may cut $X$ along $\alpha$ to get a hyperbolic sphere with two boundary components of length $\ell=\ell_\alpha(X)$, and one cusp. It's helpful to view this as a degenerate pants, where one of the boundary curves has been replaced with a cusp. 
\begin{figure}[ht!]
\includegraphics[width=.5\linewidth]{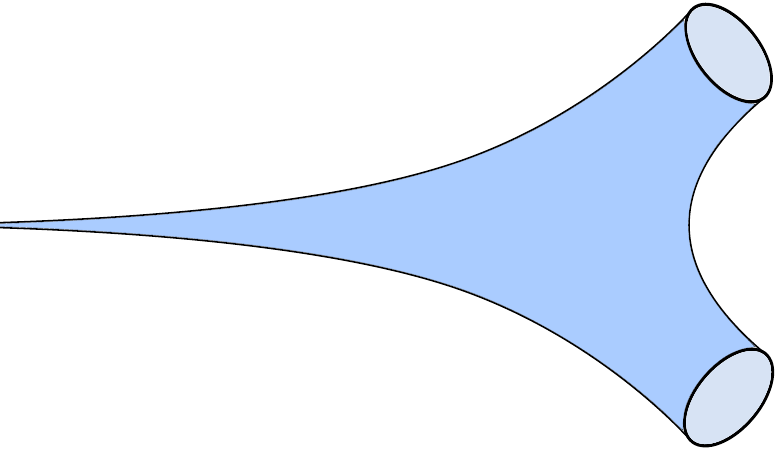}
\caption{A degenerate pants.}
\label{F:Degenerate}
\end{figure}
For any $\ell>0$, there is a unique such pants, so the point $(X,\alpha)\in \cM_{1,1}^*$ is uniquely determined by $\ell=\ell_\alpha(X)$ and a twist parameter $\tau\in [0, \ell)$ which controls how the two length $\ell$ boundary curves of the degenerate pants are glued together to give $X$. 

Using this  parametrization 
$$\cM_{1,1}^* \simeq \{(\ell, \tau): \ell>0, 0\leq \tau < \ell\}$$
and Wolpert's formula $\omega_{WP}= d \ell \wedge d\tau$, Mirzakhani concludes 
\begin{eqnarray*}
\frac12 \Vol{\cM_{1,1}} &=& \int_0^\infty \int_0^\ell \frac{1}{1+e^{\ell}} \wrt\tau \wrt\ell 
= \frac{\pi^2}{12}.
\end{eqnarray*}

\section{Integrating geometric functions over moduli space}\label{S:Integrating}

In this section we give a key result from \cite{Mirzakhani:Invent} that gives a procedure for integrating certain  functions over moduli space, generalizing the ``unfolding" step in the previous section. 

\subsection{A special case.} Let $\gamma$ be a simple non-separating closed curve on a surface of genus $g>2$. For a continuous function $f: \bR_+\to \bR_+$, we define a function $f_\gamma: \cM_{g}\to \bR$ by  
$$f_\gamma(X) = \sum_{[\alpha] \in \MCG \cdot [\gamma]} f(\ell_\alpha(X)).$$
Here the sum is over the mapping class group orbit of the homotopy class $[\gamma]$ of $\gamma$. Soon we will generalize this notation, but for the moment the subscript may seem strange: since there is only one mapping class group orbit of non-separating simple closed curve, for the moment $f_\gamma$ does not depend on $\gamma$. 

Recall that $\cM_{g-1,2}(\ell, \ell)$ is the moduli space of genus $g-1$ hyperbolic surfaces with 2 labeled boundary geodesics of length $\ell$. We will  give an outline of Mirzakhani's proof that 
$$\int_{\cM_{g}} f_\gamma(X) \dVol_{WP} = \frac12 \int_0^\infty \ell f(\ell) \Vol(\cM_{g-1,2}(\ell, \ell)) \wrt\ell.$$

Define $\cM_g^\gamma$ to be the set of pairs $(X, \alpha)$, where $X\in \cM_g$ and $\alpha$ is a geodesic with $[\alpha] \in \MCG \cdot [\gamma]$. The fibers of the map
$$\cM_g^\gamma \to \cM_g, \quad\quad (X,\alpha)\mapsto X$$
correspond exactly to the set $\{[\alpha] \in \MCG \cdot [\gamma]\}$ that is summed over in the definition of $f_\gamma$, and in fact
$$\int_{\cM_{g}} f_\gamma(X) \dVol_{WP} = \int_{\cM_g^\gamma} f(\ell_\alpha(X)) \dVol_{WP}.$$

Cutting $X$ along $\alpha$ almost determines a point of $\cM_{g-1,2}(\ell, \ell)$, except that the two boundary geodesics are not labeled. However, since there are two choices of labeling, we can say that there is a two-to-one map 
 $$\{(\ell, Y, \tau): \ell>0, Y\in\cM_{g-1,2}(\ell, \ell),  \tau\in \bR/\ell\bZ\} \to \cM_g^\gamma,$$
where the map glues together the two boundary components of $Y$ with a twist determined by $\tau$. Wolpert's Magic Formula determines the pullback of the Weil-Petersson measure, and  we get  
\begin{eqnarray*}
\int_{\cM_g^\gamma} f(\ell_\alpha(X)) \dVol_{WP}
&=& \frac12 \int_{\ell=0}^\infty \int_{\tau=0}^\ell  \int_{\cM_{g-1,2}(\ell, \ell)} f(\ell) \dVol_{WP} \wrt\tau \wrt\ell
\\&=& \frac12 \int_0^\infty \ell f(\ell) \Vol(\cM_{g-1,2}(\ell, \ell)) \wrt\ell.
\end{eqnarray*}

The case of $g=2$ is special, because  every $Y\in\cM_{1,2}(\ell,\ell)$ has an involution  exchanging the two boundary components. Because of this involution, one cannot distinguish between the two choices of labeling the two boundary components, and the map that was two-to-one is now one-to-one.  Thus, the same formula holds in genus $2$ with the factor of $\frac12$ removed. 


\subsection{The general case.}  A simple multi-curve, often just multi-curve for short, is a finite sum of disjoint simple closed curves with positive real weights, none of whose components are loops around a cusp. If $\gamma=\sum_{i=1}^k c_i \gamma_i$ is a multi-curve, its length is defined by 
$$\ell_\gamma(X) = \sum_{i=1}^k c_i \ell_{\gamma_i}(X).$$ 
We define $f_\gamma$ for multi-curves in the same way as above, and note that $f_\gamma$ in fact only depends on the mapping class group orbit of $[\gamma]$.

Suppose that cutting the geodesic representative of $\gamma$ decomposes $X\in \cM_{g,n}(L)$ into $s$ connected components $X_1, \ldots, X_s$, and that 
\begin{itemize}
\item $X_j$ has genus $g_j$,
\item $X_j$ has $n_j$ boundary components, and 
\item the lengths of the boundary components of $X_j$ are given by $\Lambda_j\in \bR_+^{n_j}$.
\end{itemize}
If we set $\ell_i = \ell_{\gamma_i}(X)$, then all the entries of each $\Lambda_j$ are from $\{L_1, \ldots, L_n\}$ (if they correspond to the original boundary of $X$) or $\{\ell_1, \ldots, \ell_k\}$ (if they correspond to the new boundary created by cutting $\gamma$).

\begin{thm}[Mirzakhani's Integration Formula]
For any multi-curve $\gamma=\sum_{i=1}^k c_i \gamma_i$, 
\begin{eqnarray*}
&&\int_{\cM_{g,n}(L)} f_\gamma \dVol_{WP}
\\&=& \IC_\gamma \int_{\ell=(\ell_1, \ldots, \ell_k)\in \bR^k_+} \ell_1 \cdots \ell_k f(c_1\ell_1+\cdots + c_k\ell_k) \prod_{j=1}^s \Vol(\cM_{g_j, n_j}(\Lambda_j)) \wrt\ell,\end{eqnarray*}
where $\IC_\gamma\in\bQ_+$ is an explicit constant.\footnote{
Slightly different values of the constant have been recorded in different places in the literature. We believe the correct constant is 
$$\IC_\gamma = \frac1{2^M  [\Stab(\gamma):\langle S, \cap_{i=1}^k \Stab^+(\gamma_i) \rangle]},$$
where $M$ is the number of $i$ such that $\gamma_i$ bounds a torus with no other boundary components and not containing any other component of $\gamma$, $\Stab(\gamma)$ is the stabilizer of the \emph{weighted} multi-curve $\gamma$,  $ \Stab^+(\gamma_i)$ is the subgroup of the mapping class group that fixes $\gamma_i$ and its orientation, and $S$ is the kernel of the action of the mapping class group on Teichm\"uller space. (Note that $S$ is trivial except in the case when $(g,n)$ is $(1,1)$ or $(2,0)$, in which case it has size two and is central.) Given two subgroups $H_1, H_2$, we write $\langle H_1, H_2\rangle$ for the subgroup they generate. 

The number $\IC_\gamma^{-1}$ arises as the degree of a measurable map from  
$$\cP= \Big\{(\ell, Y, \tau) : \ell\in \bR_+^k, Y\in \prod_{j=1}^s \cM_{g_j, n_j}(\Lambda_j), \tau\in \prod_{i=1}^k \bR/\ell_i\bZ\Big\}$$  to the space  $\cM_{g,n}^\gamma(L)$ of pairs $(X, \alpha)$, where $X\in \cM_{g,n}(L)$ and $\alpha$ is a multi-geodesic with $[\alpha] \in \MCG \cdot [\gamma]$.
 This natural map factors through the space $\cM_{g,n}^{\gamma,+}(L)$ of (ordered) tuples $(X, \alpha_1, \ldots, \alpha_k)$, where $X\in \cM_{g,n}(L)$, and the $\alpha_i$ are disjoint  oriented geodesics with $[\sum c_i \alpha_i]\in \MCG \cdot [\gamma]$. The degree of $\cP\to \cM_{g,n}^{\gamma,+}(L)$ is $2^M$ and the degree of $\cM_{g,n}^{\gamma,+}(L)\to \cM_{g,n}^{\gamma}(L)$ is the remaining index factor. 
}
\end{thm}

%

\section{Generalizing McShane's identity}\label{S:McShane}

The starting point for Mirzakhani's volume computations is the following result proven in \cite{Mirzakhani:Invent}. It relies on two explicit functions $\cD, \cR: \bR_+^3 \to \bR_+$ whose exact definitions are omitted here. 

\begin{thm}\label{T:McShane}
For any hyperbolic surface $X$ with $n$ geodesic boundary circles $\beta_1, \ldots, \beta_n$ of lengths $L_1, \ldots, L_n$, 
$$\sum_{\gamma_1, \gamma_2}  \cD(L_1, \ell_X(\gamma_1), \ell_X(\gamma_2))+
\sum_{i=2}^n \sum_{\gamma}  \cR(L_1, L_i, \ell_X(\gamma))  = L_1,
$$
where the first sum is over all pairs of  closed geodesics $\gamma_1, \gamma_2$ bounding a pants with $\beta_1$, and the second sum is over all simple closed geodesics $\gamma$ bounding a  pants with $\beta_1$ and $\beta_i$. 
\end{thm}

\begin{figure}[ht!]
\includegraphics[width=\linewidth]{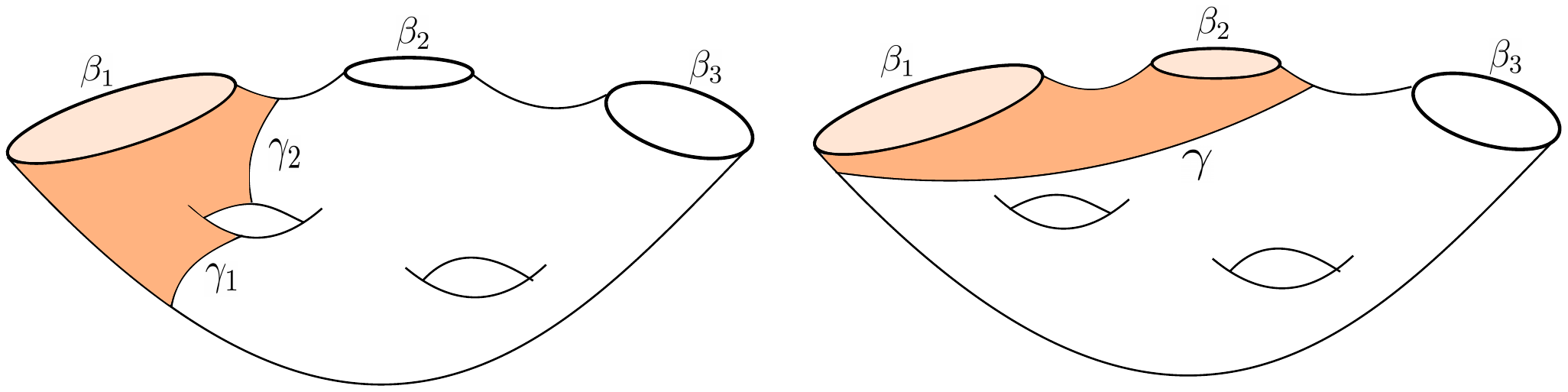}
\caption{The two types of pants in Theorem \ref{T:McShane}.}
\label{F:McPants}
\end{figure}

By studying the asymptotics of this formula when some $L_i\to 0$, it is possible to derive a related formula in the case when the boundary $\beta_i$ is replaced with a cusp. In the case when all $\beta_i$ are replaced with cusps, Mirzakhani recovers identities due to McShane \cite{McShane:Identity}, including the identity given in Section \ref{S:VolM11}. Thus, Mirzakhani refers to Theorem \ref{T:McShane} as the generalized McShane identities. 

\begin{proof}[Idea of the proof of Theorem \ref{T:McShane}]
Let $F$ be the set of points $x$ on $\beta_1$ from which the unique geodesic ray $\gamma_x$ beginning at $x$ and perpendicular to the boundary  continues forever without intersecting itself or hitting the boundary. By a result of Birman and Series, $F$ has measure 0, reflecting the fact that most geodesic rays intersect themselves \cite{BirmanSeries:Sparse}. 

It is easy to see that $\beta_1 \setminus F$ is open, and hence is a countable union of disjoint intervals $(a_h, b_h)$. 

\begin{figure}[ht!]
\includegraphics[width=0.35\linewidth]{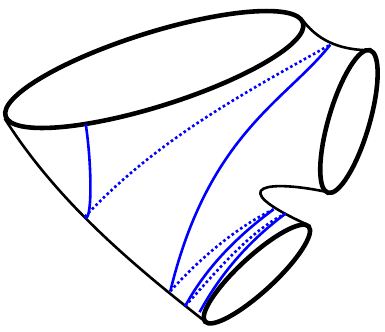}
\caption{A spiraling geodesic.}
\label{F:Spiral}
\end{figure}

Mirzakhani shows that the geodesics $\gamma_{a_h}$ and $\gamma_{b_h}$ both spiral towards either a simple closed curve or a boundary component other than $\beta_1$. There is a unique  pants $P$ with geodesic boundary containing $\gamma_{a_h}$ and $\gamma_{b_h}$. 

Each pants $P$ is associated with one or more intervals $(a_h, b_h)$, and the sum of the lengths of these intervals characterizes the functions $\cD$ and $\cR$. Having computed these functions, the identity is equivalent to $\sum_h |b_h-a_h| = L_1$. 
\end{proof}

\subsection{References.} See \cite{BridgemanTan:Identity} for a survey of related identities that have been proven since Mirzakhani's work. Of special note is that there is a related identity for closed surfaces \cite{LuoTan:Dilog}.

\section{Computation of volumes using McShane identities}\label{S:Recursive}

We now outline how Mirzakhani used her integration formula and the generalized McShane identities to recursively compute the Weil-Petersson volumes of $\cM_{g,n}(L)$ \cite{Mirzakhani:Invent}. Except in the case of $L=(0, \ldots, 0)$, $\cM_{0,4}(L)$ and $\cM_1(L)$, these volumes were unknown before Mirzakhani's work. 

As in Section \ref{S:VolM11}, we begin by integrating the generalized McShane identity to obtain 
\begin{eqnarray*}
L_1 \Vol(\cM_{g,n}(L))&=&
\int_{\cM_{g,n}(L)}\sum_{\gamma_1, \gamma_2}  \cD(L_1, \ell_X(\gamma_1), \ell_X(\gamma_2))\dVol_{WP}
\\&+& \int_{\cM_{g,n}(L)}\sum_{i=2}^n \sum_{\gamma}  \cR(L_1, L_i, \ell_X(\gamma))\dVol_{WP}.
\end{eqnarray*}
In fact, $\frac{\partial}{\partial x} \cD(x,y,z)$ and $\frac{\partial}{\partial x} \cR(x,y,z)$ are nicer functions than $\cD$ and $\cR$, so Mirzakhani considers the $\frac{\partial}{\partial L_1}$ derivative of this identity. 

Let us consider just the sum 
$$\sum_{\gamma}  \frac{\partial}{\partial L_1}\cR(L_1, L_2, \ell_X(\gamma))$$
over all simple closed geodesics $\gamma$ which bound a pants with $\beta_1$ and $\beta_2$. The set of such $\gamma$ is one  mapping class group orbit, so we may apply Mirzakhani's Integration Formula to get  
\begin{eqnarray*}
&&\int_{\cM_{g,n}(L)} \sum_{\gamma}  \frac{\partial}{\partial L_1}\cR(L_1, L_i, \ell_X(\gamma))\dVol_{WP}
\\&=&\IC_\gamma\int_{\bR_+} \ell \frac{\partial}{\partial L_1}\cR(L_1, L_i, \ell) \Vol(\cM_{g,n-1}(\ell, L_2, \ldots, L_n)) \wrt\ell.
\end{eqnarray*}

Note that the surfaces in $\cM_{g,n-1}(\ell, L_2, \ldots, L_n)$ are smaller than those in $\cM_{g,n}(L)$ in that they have one less pants in a pants decomposition. 

The sum over pairs $\gamma_1, \gamma_2$ is similar, but more complicated because the set of  multi-curves $\gamma_1+\gamma_2$ that arise  consists of a finite but possibly large number of mapping class group orbits. 

This produces an expression for $\frac{\partial}{\partial L_1} L_1 \Vol(\cM_{g,n}(L))$ as a finite sum of integrals involving volumes of smaller moduli spaces. Mirzakhani was able to compute these integrals, allowing her to compute $\Vol(\cM_{g,n}(L))$ recursively. These computations imply in particular  that $\Vol(\cM_{g,n}(L))$ is a polynomial in the $L_i^2$ whose coefficients are positive rational multiples of powers of $\pi$, which we will reprove from a different point of view in the next section. 

\subsection{References.} Mirzakhani's recursions are concisely presented in terms of the coefficients of the polynomials in \cite[Section 3.1]{Mirzakhani:Growth} and \cite[Section 2]{MirzakhaniZograf:LargeGenus}\footnote{Due to a different convention, in these papers the volumes of $\cM_{2}$ and $\cM_{1,1}(L)$ are half what our conventions give. In the second line of \cite[Section 3.1]{Mirzakhani:Growth}, there is a typo that should be corrected as $d_0=3g-3+n-|d|$. In \cite[Equation 2.13]{MirzakhaniZograf:LargeGenus}, there is a typo that should be corrected as $| I\sqcup J|=\{2, \ldots, n\}$.}. Using these recursions to compute the volume polynomials is rather slow, because of the combinatorial explosion in high genus of the number of different moduli spaces that arise from recursively cutting along geodesics. Zograf has given a faster algorithm \cite{Zograf:Large}. 

Mirzakhani's results don't directly allow for the computation of $\Vol(\cM_{g})$. However these volumes were previously known via intersection theory. They can also be recovered via the remarkable formula 
$$2\pi i (2g-2)\Vol(\cM_{g}) = \frac{\partial\Vol(\partial\cM_{g,1})}{\partial L}(2\pi i)$$
proven in \cite{DoNorbury:Cone}. 

It would  be interesting to recompute $\Vol(\cM_{g})$ using  Mirzakhani's strategy and the identity for closed surfaces in \cite{BridgemanTan:Identity}.

Mirzakhani's recursions fit into the framework of ``topological recursions" \cite{Eynard:Short}. 


\section{Computation of volumes using symplectic reduction}\label{S:VolSymp}

We now give Mirzakhani's second point of view on Weil-Petersson volumes, from \cite{Mirzakhani:Intersection}. 

\subsection{A larger moduli space} Consider the moduli space $\widehat{\cM}_{g,n}$ of genus $g$ Riemann surfaces with $n$ geodesic boundary circles with a  marked point on each boundary circle. This moduli space  has dimension $2n$ greater than that of $\cM_{g,n}(L)$, because the length of each boundary circle can vary, and the marked point on each boundary circle can vary. 

$\widehat{\cM}_{g,n}$ admits a version of local Fenchel-Nielsen coordinates, where in addition to the usual Fenchel-Nielsen coordinates there is a length parameter $\ell_i$ for each boundary circle and a parameter that keeps track of the position the each marked point on each boundary circle. The parameters keeping track of the marked points are thought of as twist parameters. The space $\widehat{\cM}_{g,n}$ also has a Weil-Peterson form $\widehat{\omega}_{WP}$, which is still described by Wolpert's Magic Formula, meaning that it is the standard symplectic form in any system of local Fenchel-Nielsen coordinates. 

Consider now the function $\mu: \widehat{\cM}_{g,n} \to \bR_+^n$ defined by
$$\mu = \left(\frac{\ell_1^2}2, \ldots, \frac{\ell_n^2}2\right).$$
The reason for this definition will become clear in the next subsection. 

Let $S^1=\bR/\bZ$, and consider the $(S^1)^n$ action that  moves the marked points along the boundary circles.  
Each level set $\mu^{-1}(L_1^2/2, \ldots, L_n^2/2)$ is invariant under the $(S^1)^n$ action, and the quotient is the space $\cM_{g,n}(L_1, \ldots, L_n)$ with fixed boundary lengths and no marked points on the boundary. That is, 
$$\cM_{g,n}(L_1, \ldots, L_n) = \mu^{-1}(L_1^2/2, \ldots, L_n^2/2)/ (S^1)^n.$$

\subsection{Symplectic reduction} We now review a version of the Duistermaat-Heckman Theorem in symplectic geometry, as it applies to $\widehat{\cM}_{g,n}$.

The symplectic form $\widehat{\omega}_{WP}$ on $\widehat{\cM}_{g,n}$  provides a non-degenerate bi-linear form on each tangent space to $\widehat{\cM}_{g,n}$. This gives an identification between the tangent space and its dual, and hence between vector fields and one-forms. 

Any function $H$ on $\widehat{\cM}_{g,n}$ determines a one-form $dH$ and hence also a vector field $V_H$ defined via this duality. This duality is recorded symbolically as $$dH = \widehat{\omega}_{WP}(V_H, \cdot).$$  The flow in the vector field $V_H$ is called the Hamiltonian flow of $H$, and $H$ is called the Hamiltonian function. 

To begin, take $n=1$. Let $\ell_1: \widehat{\cM}_{g,1}\to \bR_+$ denote the length of the unique boundary circle, and let $\tau_1$ denote the twist coordinate giving the position of the marked point. The $S^1$ action on $\widehat{\cM}_{g,n}$ discussed above is simply given by $\tau_1\mapsto \tau_1+t\ell_1$, where $t\in \bR/\bZ$, and hence is generated by the vector field $\ell_1 \partial_{\tau_1}$. Wolpert's Magic Formula gives 
$$\widehat\omega_{WP}(\ell_1 \partial_{\tau_1}, \cdot) = \ell_1 d\ell_1.$$
If $H=\ell_1^2/2$ then $dH = \ell_1 d\ell_1$, so the $S^1$ action is Hamiltonian with Hamiltonian function $H$.

Now take $n>1$. Then the $i$-th coordinate $S^1$ action on $\widehat{\cM}_{g,n}$, which moves the position of the marked point on the $i$-th boundary circle, is Hamiltonian with Hamiltonian function $\ell_i^2/2$. Using a natural way to combine different Hamiltonian functions into a single function called the moment map, one says that the $(S^1)^n$ action on $\widehat{\cM}_{g,n}$ is Hamiltonian with moment map $\mu$ given above. 

For every $\xi =(L_1^2/2, \ldots, L_n^2/2)$, the space $\mu^{-1}(\xi)$ is the manifold parameterizing surfaces in  $\cM_{g,n}(L_1, \ldots, L_n)$ with a marked point on each boundary circle, and as we've discussed 
$$\mu^{-1}(\xi)/(S_1)^n = \cM_{g,n}(L_1, \ldots, L_n).$$
The fact that this quotient is a symplectic manifold is an instance of a general phenomenon called symplectic reduction. 

Return to the case $n=1$. Then $\mu^{-1}(\xi)$ is a principal circle bundle over $\cM_{g,1}(L_1)$. (A principal $S^1$ bundle is a bundle with an action of $S^1$ that is simply transitive on fibers.) Let's call this circle bundle $\cC_1$.  

Generalizing this to $n>1$, we see that $\mu^{-1}(\xi)\to \cM_{g,1}(L_1, \ldots, L_n)$ is a product of $n$ circle bundles $\cC_i, i=1,\ldots,n$. Here $\cC_i$ can be defined as the spaces of  of surfaces in $\cM_{g,n}(L_1, \ldots, L_n)$ together with just a single marked point on the $i$-th boundary circle (and no marked points on any of the other boundary circles).

\begin{thm}[Duistermaat-Heckman Theorem]\label{T:DH}
For any fixed $\xi$ and for $t=(t_1, \ldots, t_n)\in \bR^n$  small enough, there exists a diffeomorphism   
$$\phi_t: \mu^{-1}(\xi)/(S^1)^n \to \mu^{-1}(\xi+t)/(S^1)^n$$  
such that 
$$\phi_t^*(\omega_{WP})=\omega_{WP}+\sum_{i=1}^n t_i c_1(\cC_i),$$
where $c_1(\cC_i)$ is the first Chern class of the circle bundle $\cC_i$ over $\mu^{-1}(\xi+t)/(S^1)^n$.  Here, on the left hand side $\omega_{WP}$ refers to the Weil-Petersson form on $\mu^{-1}(\xi+t)/(S^1)^n$, and on the right hand side it refers to the Weil-Petersson form on $\mu^{-1}(\xi)/(S^1)^n$. 
\end{thm}

The reader unfamiliar with Chern classes may in fact take this theorem to be the definition for this survey; we will not use any other properties of Chern classes. 

Part of Theorem \ref{T:DH} is powered by a relative of the Darboux Theorem. The Darboux Theorem states that a neighborhood of any point in a symplectic manifold is symplectomorphic to the simplest thing you could guess it to be, namely a neighbourhood in a vector space with the standard symplectic form. 

Here, a neighborhood of $\mu^{-1}(\xi)$ is topologically $\mu^{-1}(\xi)\times (-\delta, \delta)^n$, and one can create a guess for what the symplectic form $\hat\omega_{WP}$ might look like on $\mu^{-1}(\xi)\times (-\delta, \delta)^n$, using $\omega_{WP}$ on $\mu^{-1}(\xi)/(S^1)^n$ and curvature forms for the circle bundles. The Equivariant Coisotropic Reduction Theorem, which is the relative of the Darboux Theorem we referred to, says that this guess is in fact symplectomorphic to a neighborhood of  $\mu^{-1}(\xi)$ in $\widehat{\cM}_{g,n}$.

\subsection{Computations of volumes.}
From Theorem \ref{T:DH},  $\Vol(\cM_{g,n}(L_1, \ldots, L_n))$ is a polynomial in a small neighborhood of any $(L_1, \ldots, L_n)$, and hence is globally a polynomial. By considering $\xi=(\e^2/2, \ldots, \e^2/2)$, we get 
\begin{eqnarray*} 
&&\Vol(\cM_{g,n}(L_1, \ldots, L_n))
\\&=& \frac1{(3g-3+n)!} \int_{\cM_{g,n}(L_1, \ldots, L_n)} \omega_{WP}^{3g-3+n}
\\&=&\frac1{(3g-3+n)!} \int_{\cM_{g,n}(\e, \ldots, \e)} \left(\omega_{WP}+\sum \frac{L_i^2-\e^2}{2}c_1(\cC_i)\right)^{3g-3+n}.
\end{eqnarray*}
Note that  Theorem \ref{T:DH} only directly gives this for $L_i$ close to $\e$, but since the volume is a polynomial it must in fact be true for all $L_i$. 
Note also that  $\omega_{WP}$ denotes the Weil-Petersson symplectic form on $\cM_{g,n}(L_1, \ldots, L_n)$ in the first integral and on $\cM_{g,n}(\e, \ldots, \e)$ in the second integral. 

Taking a limit as $\e\to 0$, Mirzakhani obtains 
$$\Vol(\cM_{g,n}(L_1, \ldots, L_n)) = \frac1{(3g-3+n)!} \int_{\cM_{g,n}} \left(\omega_{WP}+\sum \frac{L_i^2}{2}c_1(\cC_i)\right)^{3g-3+n}.$$
Here $\cC_i$ can be defined as the circle of points on a horocycle of size 1 about the $i$-th cusp. 

One can first interpret this integral in terms of the differential forms representing $c_1(\cC_i)$ produced by the proof of Theorem \ref{T:DH}. These differential forms, and the circle bundles $\cC_i$, extend continuously to a natural compactification $\overline{\cM}_{g,n}$ constructed by Deligne and Mumford, and so one can and typically does replace $\cM_{g,n}$ with $\overline{\cM}_{g,n}$ as the space to be integrated over. This allows a more topological interpretation  of the integral as the pairing of a class in $H^{6g-6+2n}(\overline{\cM}_g)$ with the fundamental class of $\overline{\cM}_{g,n}$.

In summary, we have the following. 

\begin{thm}\label{T:poly}
The volume of $\cM_{g,n}(L_1, \ldots, L_n)$ is a polynomial 
$$\sum_{|\alpha|\leq 3g-3+n} C_g(\alpha) L^{2\alpha}$$
whose coefficients $C_g(\alpha)$ are rational multiples of integrals  of powers of the Chern classes $c_1(\cC_i)$ and the Weil-Petersson symplectic form. Here $\alpha=(\alpha_1, \ldots, \alpha_n)$,  $|\alpha|=\sum \alpha_i$ and $L^{2\alpha} = \prod L_i^{2\alpha_i}$. 
\end{thm}

We will discuss a number of interesting results and open problems about these polynomials in Section \ref{S:Random}. 

\subsection{References} For more on the material relating to symplectic reduction, see, for example, \cite[Chapter 22, 23, 30.2]{daSilva:Lectures}. 

The work of Mirzakhani suggests some similarities between moduli spaces of Riemann surfaces and spaces of representations of surface groups into compact Lie groups modulo conjugacy. These spaces of representations are also known as character varieties or moduli spaces of stable  bundles. Mirzakhani points out the connection between her techniques and those used previously by Witten and others in this context \cite{Witten:Revisited}. See the citations in Mirzakhani's papers and the survey \cite{Jeffrey:Flat} for more details.

\section{Witten's conjecture}\label{S:Witten}

We now describe Mirzakhani's proof of Witten's conjecture \cite{Mirzakhani:Intersection}. This brings us to the algebro-geometric perspective on the coefficients $C_g(\alpha)$ from Theorem \ref{T:poly}. 

\subsection{Intersection theory.} Let $\overline{\cC_i}$ and $\overline{\omega}_{WP}$ denote the extensions of $\cC_i$ and $\omega$ from $\cM_{g,n}$ to $\overline{\cM}_{g,n}$. 

 The class $c_1(\overline{\cC_i}) \in H^2(\overline{\cM}_{g,n}, \bQ)$ is typically denoted $\psi_i$, and is much studied. One often defines $\psi_i$ as the first Chern class of  a line bundle $\cL_i$ called the relative cotangent bundle at the $i$-the marked point. 

 Wolpert showed that the cohomology class $[\overline{\omega}_{WP}]\in H^2(\overline{\cM}_{g,n})$ is equal to  $2\pi^2 \kappa_1$, where $\kappa_1 \in  H^2(\overline{\cM}_{g,n}, \bQ)$ is the much studied first kappa class \cite{Wolpert:Homology}. As a result, all of the coefficients $C_g(\alpha)$ from Theorem \ref{T:poly} are in $\bQ[\pi^2]$.  

The compactification $\overline{\cM}_{g,n}$ is an algebraic variety and a smooth orbifold, and the classes $\psi_i$ and $\kappa_1$ can be thought of as dual to (equivalence classes of) divisors, which are linear combinations of subvarieties of complex codimension 1. The intersection of two such classes, if transverse, has complex codimension 2, and similarly the intersection $\dim_\bC \overline{\cM}_{g,n}=3g-3+n$ of them, if transverse, is a finite collection of points. Integrals of a product of $3g-3+n$ of the $\psi_i$ and $\kappa_1$ classes, which up to factors are exactly the coefficients $C_g(\alpha)$, count the number of points of intersection. Thus, they are called intersection numbers. See the book \cite[Chapter 4.6]{LandoZvonkin:Graphs} for some example computations using this point of view. 

It is hard for the uninitiated to fathom how much useful information such intersection numbers can contain, so we pause to give just a few points of motivation. 
\begin{itemize}
\item They are central to the study of the geometry of $\overline{\cM}_{g,n}$. 
\item By  Theorem \ref{T:DH} they determine Weil-Petersson volumes. Later we will see that these volumes can be used to understand the geometry of Weil-Petersson random surfaces. 
\item They appear in theoretical physics \cite{Witten:TwoD}. 
\item They determine counts of combinatorial objects called ribbon graphs \cite{Kontsevich:Airy}.
\item  They determine Hurwitz numbers, which count certain branched coverings of the sphere, or equivalently factorizations of permutations into transpositions \cite{EkedahlLandoShapiroVainshtein:Hurwtiz}.
\end{itemize}

\subsection{A generating function for intersection numbers.} Make the notational convention 
$$\langle \tau_{d_1}\cdots \tau_{d_n}\rangle_g = \int_{\overline{\cM}_{g,n}} \psi_1^{d_1} \ldots \psi_n^{d_n}.$$
Unless $\sum d_i=3g-3+n $, this is defined to be zero. Note that ``$\langle \tau_{d_1}\cdots \tau_{d_n}\rangle_g$" should be considered as a single mathematical symbol, and the order of the $d_i$'s doesn't matter. 

 Define the generating function for top intersection products in genus $g$ by 
$$F_g(t_0, t_1, \ldots) = \sum_n \frac1{n!}\sum_{d_1, \ldots, d_n} \langle \prod \tau_{d_i}\rangle_g t_{d_1} \cdots t_{d_n},$$
where the sum is over all non-negative sequences $(d_1, \ldots, d_n)$ such that $\sum d_i = 3g-3+n$. 
One can then form the generating function 
$$F=\sum_g \lambda^{2g-2} F_g,$$
which arises as a partition function in 2D quantum gravity. Note that $F$ is a generating function in infinitely many variables: $\lambda$ keeps track of the genus, and $t_d$ keeps track of the number of $d$-th powers of psi classes. 

Witten's conjecture is equivalent to the fact that $e^F$ is annihilated by a sequence of differential operators 
$$L_{-1}, L_0, L_1, L_2, \ldots$$
satisfying the Virasoro relations 
$$[L_m, L_k] = (m-k)L_{m+k}.$$ 
To give an idea of the complexity of these operators, we record the formula for $L_n, n>0$: 
\begin{eqnarray*}
L_n 
=&&
 -\frac{(2n+3)!!}{2^{n+1}}\frac{\partial}{\partial t_{n+1}}
 \\&+& \frac{1}{2^{n+1}}\sum_{k=0}^\infty \frac{(2k+2n+1)!!}{(2k-1)!!}t_k \frac{\partial}{\partial t_{n+k}}
 \\&+& \frac1{2^{n+2}} \sum_{i+j=n-1} (2i+1)!! (2k+1)!! \frac{\partial^2}{\partial t_{i}\partial t_{j}}.
 \end{eqnarray*}

The equations $L_i(e^F)=0$ encode recursions among the intersection numbers, which appear as the constant terms in Mirzakhani's volume polynomials. These recursions allow for the computation of all intersection numbers of psi classes. Mirzakhani showed that these recursions follow from her recursive formulas for the volume polynomials, thus giving a new proof of Witten's conjecture. 

\subsection{A brief history.} Witten's conjecture was published in 1991, motivated by physical intuition that two different models for 2D quantum gravity should be equivalent \cite{Witten:TwoD}. Kontsevich published a proof in  1992, using a combinatorial model for $\cM_{g,n}$ arising from Strebel differentials, ribbon graphs, and random matrices \cite{Kontsevich:Airy}. This work was central in his 1998 Fields Medal citation. 

 It wasn't until 2007 that Mirzakhani's proof was published, and around the same time other proofs appeared. Later, Do related Mirzakhani's and Kontsevich's proofs, recovering Kontsevich's formula for the number of ribbon graphs by considering asymptotics of the Weil-Petersson volume polynomials, using that a rescaled Riemann surface with very large geodesic boundary looks like a graph \cite{Do:Asymptoptic}.

%

\section{Counting simple closed geodesics}\label{S:CountSimple}

Let $X$ be a complete hyperbolic surface, and let $c_X(L)$ be the number of primitive closed geodesics of length at most $L$ on $X$. Primitive means simply that the geodesic does not traverse the same path multiple times. The famous Prime Number Theorem for Geodesics gives the asymptotic 
$$c_X(L) \sim \frac12 \frac{e^L}{L}$$
as $L\to \infty$. (The factor of $\frac12$ disappears if one counts primitive oriented geodesics, since there are two orientations on each closed geodesic.)  Amazingly, this doesn't depend on which surface $X$ we choose, or even the genus of $X$.

 In \cite{Mirzakhani:Annals}, Mirzakhani proved that the number of closed geodesics of length at most $L$ on $X$ that don't intersect themselves   is asymptotic to a constant depending on $X$ times $L^{6g-6+2n}$.  That this asymptotic is polynomial rather than exponential reflects the extreme unlikeliness that a random closed geodesic is simple, in the same spirit as the result of Birman and Series mentioned in Section \ref{S:McShane}. 

In fact Mirzakhani proved a more general result. For any rational multi-curve $\gamma$, she considered 
$$s_X(L,\gamma) = |\{\alpha \in \MCG \cdot \gamma : \ell_\alpha(X) \leq L\}|.$$
In other words, $s_X(L, \gamma)$ counts the number of closed multi-geodesics $\alpha$ on $X$ of length less than $L$ that are ``of the same topological type" as $\gamma$. 

The set of simple closed curves forms finitely many mapping class group orbits. So by summing finitely many of these functions $s_X(L,\gamma)$, one gets the corresponding count for all simple closed curves. 

\begin{thm}\label{T:Ann}
For any rational multi-curve $\gamma$, 
$$\lim_{L\to \infty} \frac{s_X(L, \gamma)}{L^{6g-6+2n}} = \frac{c(\gamma) \cdot B(X)}{b_{g,n}},$$
where $c(\gamma)\in \bQ_+$, $B :\cM_{g,n} \to \bR_+$ is a proper, continuous function with a simple geometric definition, and $b_{g,n}=\int_{\cM_{g,n}} B(X) \dVol_{WP}$. 
\end{thm} 


Eskin-Mirzakhani-Mohammadi have recently given a new proof of Theorem \ref{T:Ann} that gives an error term, which we will discuss in Section \ref{S:Effective}, and Erlandsson-Souto have also given a new proof \cite{ErlandssonSouto:MCC}. Here we outline the original proof, after first commenting on one application. 

\subsection{Relative frequencies.} Consider, for example, the case $(g,n)=(2,0)$ of closed genus 2 surfaces, just to be concrete. The set of simple closed curves consists of two mapping class groups orbits:  the orbit of a non-separating curve $\gamma_{ns}$ and the orbit of a curve $\gamma_{sep}$ that separates the surface into two genus one subsurfaces. 

The fact that the limit in Theorem \ref{T:Ann} is the product of a function of $\gamma$ and a function of $X$ has the following  consequence: A very long simple closed curve on $X$, chosen at random among all such curves, has probability about 
$$\frac{c(\gamma_{sep})}{c(\gamma_{ns})+c(\gamma_{sep})}$$
of being separating. Remarkably, this probability is computable and does not depend on $X$! 

Even more remarkably, recent discoveries prove that the same probabilities appear in discrete problems about surfaces assembled out of finitely many unit squares \cite{DelecroixGoujardZografZorich,Arana-Herrera:Counting}. 

\subsection{The space of measured foliations.} The space of rational multi-curves  admits a natural completion called the space $\MF$ of measured foliations. Later we will  delve into measured foliations, but here we only need a few properties of this space. 
\begin{itemize}
\item $\MF$ is homeomorphic to $\bR^{6g-6+2n}.$
\item $\MF$ does not carry a natural linear  structure. The most superficial indication of this is that any closed curve $\alpha$ gives a point of $\MF$, but there is no ``$-\alpha$" in $\MF$, because multi-curves are defined to have positive coefficients. There is however a natural action of $\bR_+$ on $\MF$, which on multi-curves simply multiplies the coefficients by $t\in \bR_+$. 
\item $\MF$ has a natural piece-wise linear integral structure, that is, an atlas of charts to $\bR^{6g-6+2n}$ whose transition functions are piece-wise in $GL(n,\bZ)$. 
\item Define the integral points of $\MF$ as the set $\MF(\bZ)\subset \MF$ of points mapping to $\bZ^{6g-6+2n}$ under the charts. Define the rational points $\MF(\bQ)$ similarly. Then integral (resp. rational) points of $\MF$ parametrize homotopy classes of integral (resp. rational) multi-curves on the surface. 
\item Any $X\in \cM_{g,n}$ defines a continuous length function $$\MF\to \bR_+, \quad\quad \lambda\to \ell_\lambda(X)$$ whose restriction to multi-curves gives the hyperbolic length of the geodesic representative of the multi-curve on $X$. In particular, $\ell_{t\lambda}(X) = t\ell_\lambda(X)$ for all $t\in \bR$. 
\end{itemize} 

\subsection{Warm up.} How many points of $\bZ^2\subset \bR^2$ have length at most $L$? It is equivalent to ask  about the number of points of $\frac1{L} \bZ^2$ contained in the unit ball. 

Recall that the Lebesgue measure  can be defined as the limit as $L\to \infty$ of 
$$\frac1{L^2} \sum_{\alpha\in \bZ^2} \delta_{\frac1{L} \alpha},$$
where $\delta_x$ denotes the point mass at $x$.
Hence the number of points of $\frac1{L} \bZ^2$ contained in the unit ball is asymptotic to $L^2$ times the Lebesgue measure of the unit ball. 

\subsection{The Thurston measure.} Let's start with the easy question of asymptotics for the number 
$$S_X(L) = |\{\alpha \in \MF(\bZ) : \ell_\alpha(X) \leq L\}|$$
of all integral multi-curves of length at most $L$. Using $\ell_{t\lambda}(X) = t\ell_\lambda(X),$ we observe that
$$S_X(L) = |\{\alpha \in L^{-1}\MF( \bZ) : \ell_\alpha(X) \leq 1\}|.$$ 
It's now useful to define the ``unit ball" 
$$B_X = \{\alpha \in \MF :  \ell_\alpha(X) \leq 1\}$$
and the measures 
$$\mu^L = \frac{1}{L^{6g-6+2n}} \sum_{\alpha\in \MF(\bZ)} \delta_{\frac1{L}\alpha}.$$
With these definitions, 
$$S_X(L)= L^{6g-6+2n} \mu^L(B_X).$$
As in our warm up, the measures $\mu^L$ converge to a natural Lebesgue class measure $\mu_{Th}$ on $\MF$. This measure $\mu_{Th}$ is called the Thurston measure and is Lebesgue measure in the charts mentioned above. If we define $B(X)= \mu_{Th}(B_X)$, we get the asymptotic 
$$S_X(L) \sim B(X) L^{6g-6+2n}.$$

\subsection{The proof.}
Mirzakhani's approach to Theorem \ref{T:Ann} similarly defines measures 
$$\mu^L_\gamma = \frac{1}{L^{6g-6+2n}} \sum_{\alpha\in \MCG \cdot \gamma} \delta_{\frac1{L}\alpha}.$$
As above, to prove Theorem \ref{T:Ann}, it suffices to show the convergence of measures  
$$\mu^L_\gamma \to \frac{c(\gamma)}{b_{g,n}} \mu_{Th}.$$
 
Using the Banach-Alaoglu Theorem, it isn't hard to show that there are subsequences $L_i\to \infty$ such that $\mu^{L_i}_\gamma$ converges to some measure $\mu^\infty_\gamma$, which might a priori depend on which subsequence we pick. To prove Theorem \ref{T:Ann} it suffices to show that, no matter which such subsequence we use, we have 
$$\mu^\infty_\gamma = \frac{c(\gamma)}{b_{g,n}} \mu_{Th}.$$ 

By definition,  $\mu^L_\gamma\leq \mu^L$, since the mapping class group orbit of $\gamma$ is a subset of $\MF(\bZ)$. Since $\mu^L$ converges to the Thurston measure, we get that $\mu^\infty_\gamma \leq  \mu_{Th}$. 

Since $\mu^L_\gamma$ is mapping class group invariant, the same is true for $\mu^\infty_\gamma$. A result of Masur in ergodic theory, which we will discuss in Section \ref{S:Ergodic}, gives that any mapping class group invariant measure on $\MF$ that is absolutely continuous to $\mu_{Th}$ must be a multiple of $\mu_{Th}$   \cite{Masur:Ergodic}.  

So $\mu^L_\gamma\leq \mu^L = c \mu_{Th}$ for some $c\geq 0$. At this point in the argument as far as we know $c$ could depend on the subsequence of $L_i$.

Unraveling the definitions, we have that  
\begin{equation}\label{E:sX}
\frac{s_X(L_i, \gamma)}{L_i^{6g-6+2n}} \to c\cdot B(X),
 \end{equation} 
for any $X\in \cM_{g,n}$. Writing 
$$s_X(L_i, \gamma)= \sum_{\alpha \in \MCG \gamma} \chi_{[0,L_i]} (\ell_\alpha(X)),$$
we recognize the type of function that Mirzakhani's Integration Formula applies to. By integrating the left hand side of \eqref{E:sX} over moduli space, Mirzakhani is able to prove that $c=\frac{c(\gamma)}{b_{g,n}}$ as desired. On the one hand, the integral of $c B(X)$ is $c\cdot b_{g,n}$. On the other hand, the limit  $c(\gamma)$ of the integral of $\frac{s_X(L_i, \gamma)}{L_i^{6g-6+2n}}$ is easily expressed in terms of the leading order term in one of Mirzakhani's volume polynomials. 

\subsection{Open problems.} We will return to counting later, but for now we mention the following. 
\begin{prob}
Prove an analogue of Theorem \ref{T:Ann} for non-orientable hyperbolic surfaces. 
\end{prob}
An example is known already with asymptotics $L^\delta$ with $\delta$ non-integral \cite{Magee:OneSided}.  See \cite{Gendulphe:Wrong} for a more precise conjecture, as well as a number of related open problems and an analogy between moduli spaces of non-oriented hyperbolic surfaces and infinite volume geometrically finite hyperbolic manifolds. 

%
%

\section{Random surfaces of large genus}\label{S:Random}

Given a random $d$-regular graph with many vertices, what is the chance that it contains a short loop? Is a random graph easy to cut in two? What properties can be expected of the graph Laplacian?  

Mirzakhani considered analogues of these well-studied questions for Weil-Petersson random Riemann surfaces \cite{Mirzakhani:Growth, MirzakhaniZograf:LargeGenus, MirzakhaniPetri:Lengths}, and devoted her 2010 talk at the International Congress of Mathematicians to this topic \cite{Mirzakhani:ICM}. 

In this section we discuss this work. We will leave out the background on graphs, but many readers will wish to keep in mind the comparison between a random $d$-regular graph, with $d$ fixed and a large number of vertices, and a random surface with large genus. 

\subsection{Understanding the volume polynomials.} We begin with the constant term of the polynomial $\Vol(\cM_{g,n}(L))$, which is  the volume $V_{g,n}$ of $\cM_{g,n}$. Improving on previous results of Mirzakhani and others, Mirzakhani and Zograf proved the following \cite{Mirzakhani:Growth, MirzakhaniZograf:LargeGenus}.

\begin{thm}\label{T:asym}
There exists a universal constant $C\in (0,\infty)$ such that for any fixed $n$, $V_{g,n}$ is asymptotic to 
$$C \frac{(2g-3+n)! (4\pi^2)^{2g-3+n}}{\sqrt{g}}$$
as $g\to\infty$.
\end{thm}
This largely verified a previous conjecture of Zograf, except that his prediction that $C=\frac1{\sqrt{\pi}}$ is still open \cite{Zograf:Large}. Mirzakhani and Zograf also gave a more detailed asymptotic expansion. The proof uses the recursions satisfied by $V_{g,n}$ discussed in Section \ref{S:Recursive}.  

Previous results gave asymptotics  as $n\to \infty$ for fixed $g$ \cite{ManinZograf:Invertible}. See \cite[Section 1.4]{Mirzakhani:Growth} for open questions concerning asymptotics as both $g$ and $n$ go to infinity. 

Also by studying recursions, Mirzakhani proved results in \cite{Mirzakhani:Growth} that imply
\begin{equation}\label{E:sinh}
\Vol(\cM_{g,n}(L)) \leq V_{g,n} \prod_{i=1}^n \frac{\sinh(L_i/2)}{L_i/2}.
\end{equation}
Mirzakhani and Petri showed this bound is asymptotically sharp for fixed $n$ and bounded $L$ as $g\to\infty$ \cite[Proposition 3.1]{MirzakhaniPetri:Lengths}. The proof of the inequality actually gives a bound with $\sinh$ replaced with one of its Taylor polynomials. 


\subsection{An example.} To illustrate Mirzakhani's techniques, we will give an upper bound for the probability that a random surface in $\cM_g$ has a non-separating simple closed geodesic of length at most some small $\e>0$.  

We begin by studying the average over $\cM_g$ of the number of simple, non-separating geodesics of length at most $\e$ on $X\in \cM_{g}$. If $\gamma$ is a simple non-separating curve, we can express this as 
$$\frac{1}{V_g}\int_{\cM_g} \sum_{\alpha\in \MCG\cdot\gamma} \chi_{[0,\e]} (\ell_\alpha(X)) \dVol_{WP},$$
where $\chi_{[0,\e]}$ is the characteristic function of the interval $[0,\e]$. Mirzakhani's Integration Formula gives that this is equal to a constant times
$$ \frac{1}{V_g}\int_0^\e \ell \Vol(\cM_{g-1,2}(\ell, \ell)) \wrt\ell.$$
Since $\ell$ is small, inequality \eqref{E:sinh} gives that $\Vol(\cM_{g-1,2}(\ell, \ell))$ is approximately equal to the constant term $V_{g-1,2}$ of the volume polynomial, so the average is approximately a constant times
$$\frac{V_{g-1,2}}{V_g} \e^2.$$
The asymptotics in Theorem \ref{T:asym} imply that $\frac{V_{g-1,2}}{V_g}$ converges to 1 as $g\to\infty$, so we get that the average number of simple, non-separating geodesics of length at most $\e$ is asymptotic, as $g\to\infty$, to a constant times $\e^2$. In particular, this implies that the probability that a random surface in $\cM_g$ has such a geodesic is bounded above by a constant times $\e^2$. 

A similar lower bound is possible by giving upper bounds for the average number of pairs of non-separating simple closed curves. 

\subsection{Results.} Here is an overview of results from \cite{Mirzakhani:Growth}, which concern random $X\in \cM_g$ as $g\to\infty$. 
\begin{itemize}
\item The probability that $X$ has a geodesic of length at most $\e$ is bounded above and below by a constant times $\e^2$. 
\item The probability that $X$ has a separating geodesic of length at most $1.99\log(g)$  goes to 0. 
\item The probability that $X$ has Cheeger constant  less than $0.099$ goes to 0. 
\item The probability that $\lambda_1(X)$, the first eigenvalue of the Laplacian, is less than $0.002$ goes to 0. 
\item The probability that the diameter of $X$ is greater than $40\log(g)$ goes to 0.
\item The probability that $X$ has an embedded ball of radius at least $\log(g)/6$ goes to 1. 
\end{itemize}

The first two results are proven using the techniques in the example. The Cheeger constant is defined as
$$h(X) = \inf_\alpha \frac{\ell(\alpha)}{\min(\Area(X_1), \Area(X_2))},$$
where the infimum is over all smooth multi-curves $\alpha$ that cut $X$ into two subsurfaces $X_1, X_2$. Mirzakhani defines the geodesic Cheeger constant $H(X)$ to be the same quantity where $\alpha$ is required to be a geodesic multi-curve, so obviously $h(X)\leq H(X)$. She proves that 
$$\frac{H(X)}{H(X)+1}\leq h(X),$$
and is then able to study $H(X)$ using the techniques in the example. The result on $\lambda_1$ follows from the Cheeger inequality $\lambda_1\geq h(X)^2/4$. 


We conclude with a special case of the main result of Mirzakhani and Petri \cite{MirzakhaniPetri:Lengths}. 

\begin{thm}
For any $0<a<b$, the number of primitive closed geodesics of length in 
 $[a,b]$, viewed as a random variable on $\cM_g$, converges to a Poisson
distribution as $g\to\infty$. 
\end{thm}

What is fascinating about this result of Mirzakhani and Petri is that it concerns all primitive closed geodesics, not just the simple ones. The proof uses that a geodesic $\gamma$ of length at most a constant  $b$ on a surface $X$ of very large genus  is contained in a subsurface of bounded genus and with a bounded number of boundary components (depending on $b$). The boundary of that subsurface is a simple multi-curve $\beta$ associated to $\gamma$. By showing that, as $g\to\infty$, most $X$ do not have a separating multi-curve of bounded length, they are able to show that on most $X$ most primitive geodesics are simple, and hence use the techniques illustrated in the example. 

\subsection{Open problems.} For some problems, we list an easier version followed by a harder version. 

\begin{prob}
Does there exist a sequence of Riemann surfaces $X_n$ of genus going to infinity with $\lambda_1(X_n)\to \frac14$?
Does $\lambda_1$ converge to $\frac14$ in probability as $g\to \infty$? 
\end{prob}

\begin{prob}
Is it true that for all $g$ there is an $X\in \cM_g$ such that $\lambda_1(X)>1/4$? 
Is $\liminf_{g\to\infty} \Prob(\lambda_1>\frac14)>0$?  
\end{prob}

\begin{prob}
Is there an $\e>0$ so that $\Prob(h\leq1-\e)\to 1$ as $g\to \infty$? 
Is there an $\e>0$ so that there are no surfaces with $h>1-\e$ in sufficiently high genus? 
\end{prob}

Following conversations with Mike Lipnowski, the author finds it plausible that all three problems have a positive answer. A version of the first part of the first problem appears as \cite[Conjecture 5]{WuXue:Small}. 

Mirzakhani was also interested discrete models of random surfaces, resulting from gluing together triangles \cite{BrooksMakover:Random}, and in the collection of  all covers of a fixed surface. 

\begin{prob}
Fix $X\in \cM_2$ and $\e>0$. Is there a $C>0$ such that for every $g>2$, every $Y\in \cM_g$ with no geodesic of length less than $\e$ has Teichm\"uller distance at most $C$ from some (unramified) cover of $X$? 
\end{prob}

\section{Preliminaries on dynamics on moduli spaces}\label{S:PreTD}

This section will introduce the central concepts for the remainder of our tour. 

\subsection{Polygonal presentations of quadratic differentials.}\label{SS:poly} 
\begin{figure}[ht!]
\includegraphics[width=0.5\linewidth]{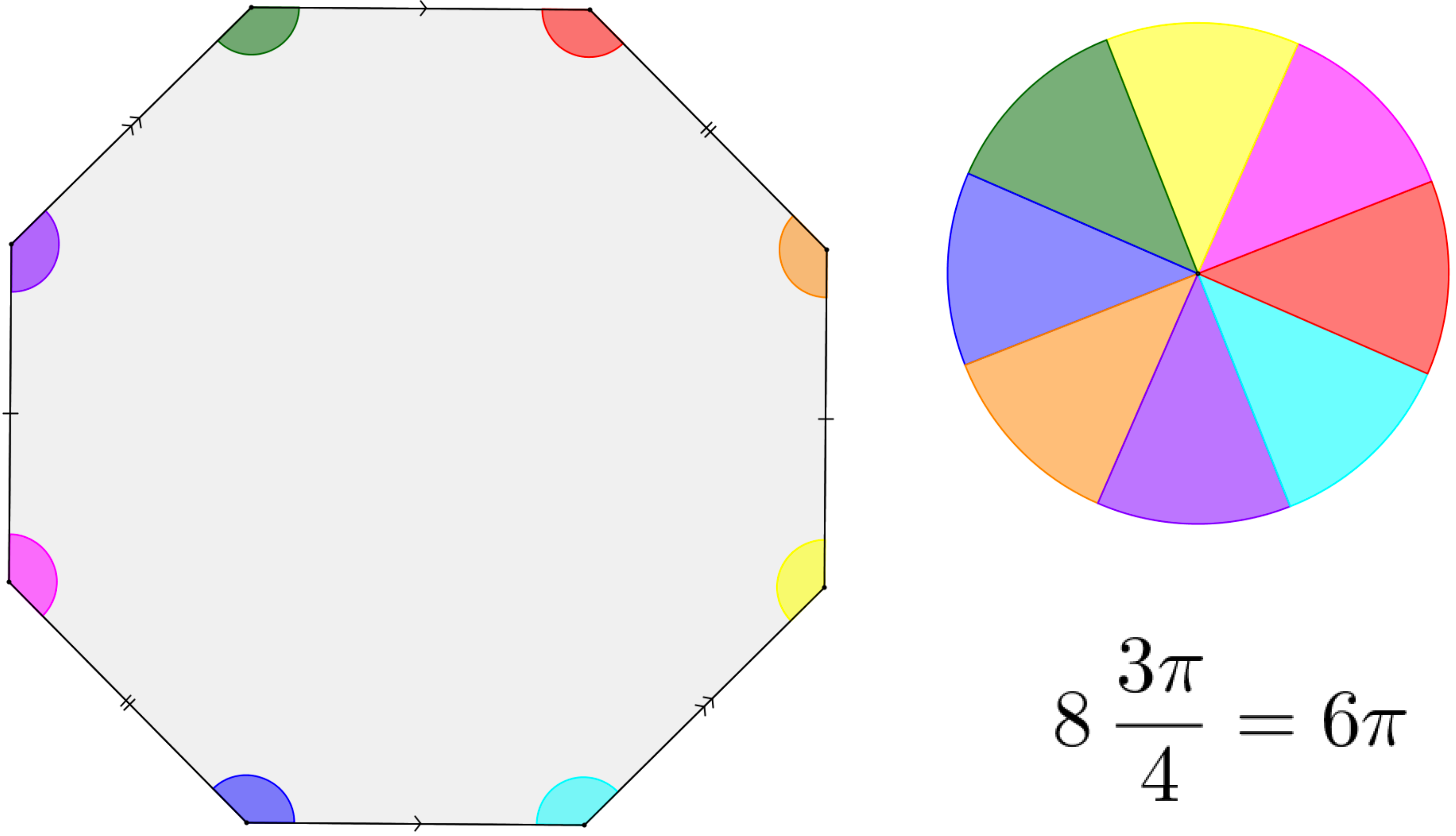}
\caption{The octagon with opposite sides identified.}
\label{F:ConeAngle}
\end{figure}

Consider the regular octagon with opposite sides identified, as in Figure \ref{F:ConeAngle}. It defines a genus 2 surface with a flat metric, except that the metric has a cone point singularity with $6\pi$ angle at the single point of the surface resulting from identifying all eight vertices. 

More generally, consider any collection of disjoint polygons in $\bC$, and glue parallel edges via maps of the form $z\mapsto \pm z + C$ to obtain a closed surface $X$. This surface is flat away from cone points created by identifying vertices. On the complement of the set $\Sigma$ of cone points, the surface has an atlas of charts to $\bC$ with transition functions of the form $z\mapsto \pm z+C$, and one can prove that the cone angle at every vertex is an integer multiple of $\pi$. Using the atlas of charts, we see that the differential $(dz)^2$ is well-defined on $X-\Sigma$, and one can further prove that it extends to a meromorphic differential with a zero of order $k$ at each cone point of angle $(2+k)\pi$. Here we count simple poles as zeros of order $k=-1$, and the differential is holomorphic away from the simple poles. 

The resulting differential is called a quadratic differential. It can be defined as a section of a complex line bundle over $X$, but the polygonal point of view will suffice for much of our discussion. Indeed, foundational results show that every non-zero quadratic differential has a polygonal presentation as above, and two polygonal presentations define the same quadratic differential if and only if they are related by a sequence of cut and paste moves. Typically it will be implicit that the quadratic differentials we discuss are non-zero. 

The simplest quadratic differential is $(dz)^2$ defined on the complex plane. Since it is invariant under translations, it descends to give a quadratic differential, which we will also call $(dz)^2$, on the torus $\bC/\bZ[i]$. A polygonal presentation is the 1 by 1 square with opposite sides identified. 

See, for example, \cite[Section 1]{Wright:Broad} for more details on the material in this subsection and the next. 

\subsection{Moduli spaces.}\label{SS:moduli} 
\begin{figure}[ht!]
\includegraphics[width=.8\linewidth]{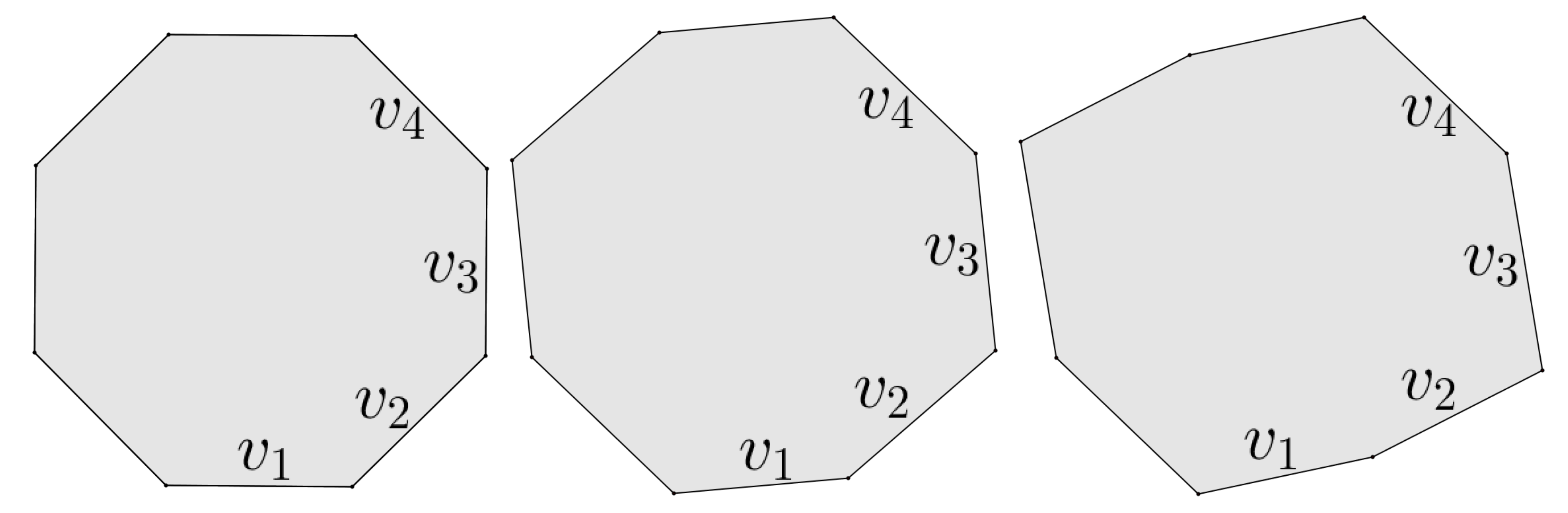}
\caption{Small deformations of the octagon surface are parametrized by $(v_1, v_2, v_3, v_4)\in \bC^4$.}
\label{F:PeriodCoords}
\end{figure}

Consider the regular octagon with opposite sides identified. We can deform this quadratic differential by changing four of the edge vectors, which each may be viewed as an element of $\bC$, as in Figure \ref{F:PeriodCoords}. We pick these four edges to contain one edge out of every pair of edges that are identified to each other. Since paired edges must have the same length and direction, the deformation is specified by the change to just the four edge vectors, and we can  guess that the deformation space is locally parametrized by a small open set in $\bC^4$. 

If a quadratic differential has $s$ zeros of order $\kappa=(\kappa_1, \ldots, \kappa_s)$, or equivalently cone points of order $((2+\kappa_1)\pi, \ldots, (2+\kappa_s)\pi)$, then the genus of the surface is given by $4g-4=\sum \kappa_i$. 

Define the stratum $\cQ(\kappa)$ as the set of all  quadratic differentials with zeros of order $\kappa$. It turns out that $\cQ(\kappa)$ is a complex orbifold, and moreover it has an atlas of charts to $\bC^m$ with transition functions in $GL(m, \bZ)$. Each  $\cQ(\kappa)$ has finitely many connected components, and they have complex dimension either $m=2g+s-1$ or $m=2g+s-2$. The $m$ local coordinates $\bC^m$  can be thought of as $m$ edge vectors in a polygonal presentation, and the change of coordinate functions in   $GL(m, \bZ)$ correspond to doing a cut and paste and picking new edge vectors. 

\begin{figure}[ht!]
\includegraphics[width=.6\linewidth]{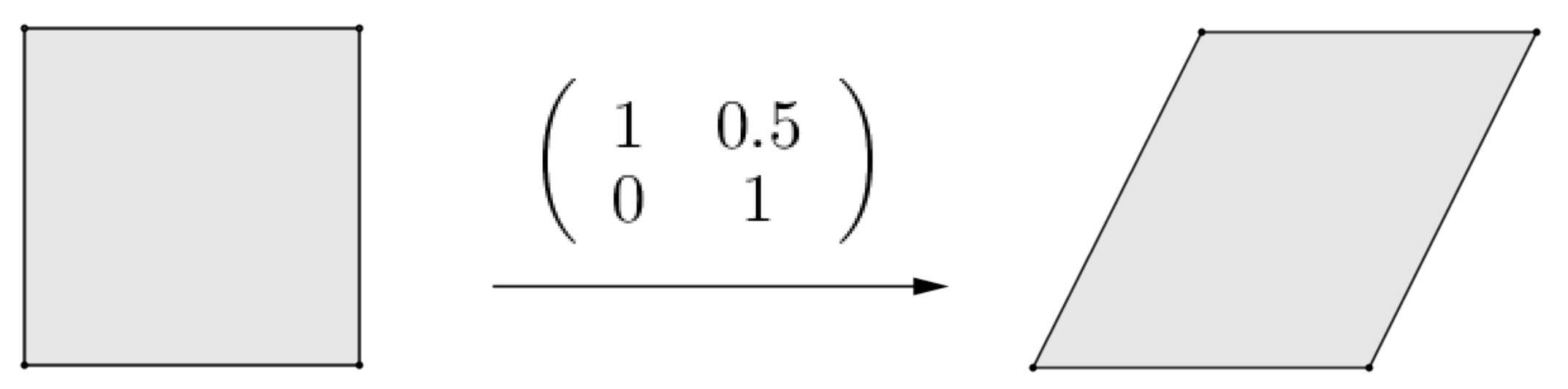}
\caption{If $g\in GL(2,\bR)$ and $(X,q)\in \cQ(\kappa)$, then $g(X,q)$ is defined by letting $g$ act on a polygonal presentation of $(X,q)$ to obtain a polygonal presentation of $g(X,q)$. In this example, $(X,q) = (\bC/\bZ[i], (dz)^2)$.}
\label{F:Shear}
\end{figure}

The $GL(2,\bR)$ action on $\bC\simeq \bR^2$ induces an action of $GL(2, \bR)$ on each stratum $\cQ(\kappa)$, as in Figure \ref{F:Shear}. One often considers just the action by the connected subgroup $\GL\subset GL(2,\bR)$ of matrices with positive determinant. 

\subsection{Abelian differentials.} This subsection is a prerequisite only for Section \ref{S:OrbitClosures}.
 
Consider a quadratic differential $(X,q)$ obtained from polygons, and suppose all the edge identifications are via maps of the form $z\mapsto  z + C$, rather than $z\mapsto  \pm z + C$. In this case, the quadratic differential $q$ is the square of an Abelian differential $\omega$, that is $q=\omega^2$. An Abelian differential is a holomorphic one-form. In the coordinates provided by the complex plane, away from the cone points, $q=(dz)^2$ and $\omega=dz$.

One can define strata $\cH(\kappa)$ of Abelian differentials as for strata of quadratic differentials, and there is a $GL(2,\bR)$  action on $\cH(\kappa)$. Up to passing to a double cover, every quadratic differential is the square of an Abelian differential, so it is often possible to study strata of Abelian differentials rather than quadratic differentials.  

Again $\cH(\kappa)$ has local coordinates given by edge vectors in a polygonal presentation of the surface, but now these coordinates also admit an additional interpretation. Namely, each edge defines a relative homology class $\gamma\in H_1(X,\Sigma)$, and the corresponding coordinate is the relative period
$$\int_\gamma \omega.$$
Using a basis $\gamma_1, \ldots, \gamma_m$ of $H_1(X,\Sigma)$ gives the local coordinates, now called period coordinates. 

The local coordinates can be thought of as the composition of the map 
$$(X,\omega)\mapsto [\omega]\in H^1(X,\Sigma, \bC)$$ 
with an isomorphism $H^1(X,\Sigma, \bC)\simeq \bC^m$. Here $[\omega]$ denotes the relative cohomology class of $\omega$. 

\subsection{Relationship to Teichm\"uller Theory.}\label{SS:Teich} 
Orbits of 
$$g_t = \left( \begin{array}{cc} e^t & 0\\ 0 & e^{-t} \end{array}\right)\subset \GL$$ 
project via $(X,q)\to X$ to geodesics in $\cM_g$ for a natural  metric on $\cM_{g,n}$ called the Teichm\"uller metric. Orbits of $\GL$ project to holomorphic and isometric immersions of the hyperbolic plane $\bH$ into $\cM_g$ called Teichm\"uller discs or complex geodesics. 

The Teichm\"uller distance $d(X,Y)$ between two Riemann surfaces $X,Y\in \cM_{g,n}$ measures how non-conformal a map $X\to Y$ must be. The Teichm\"uller metric is complete, and each pair of points in $\cT_{g,n}$ are joined by a unique Teichm\"uller geodesic. 

\subsection{Measured foliations and laminations.}\label{SS:MFML} This subsection is a prerequisite only for Sections \ref{S:Earthquake} and \ref{S:Ergodic}.

 Every quadratic differential $(X,q)$ defines a flat metric with cone singularities on the surface, but in fact it defines a bit more structure than that. Note that the atlas of charts away from the singularities have transition functions of the form  $z\mapsto \pm z+C$, and a general isometry of $\bR^2$ does not have this form. For example, rotations are not allowed as transition functions. 

\begin{figure}[ht!]
\includegraphics[width=.6\linewidth]{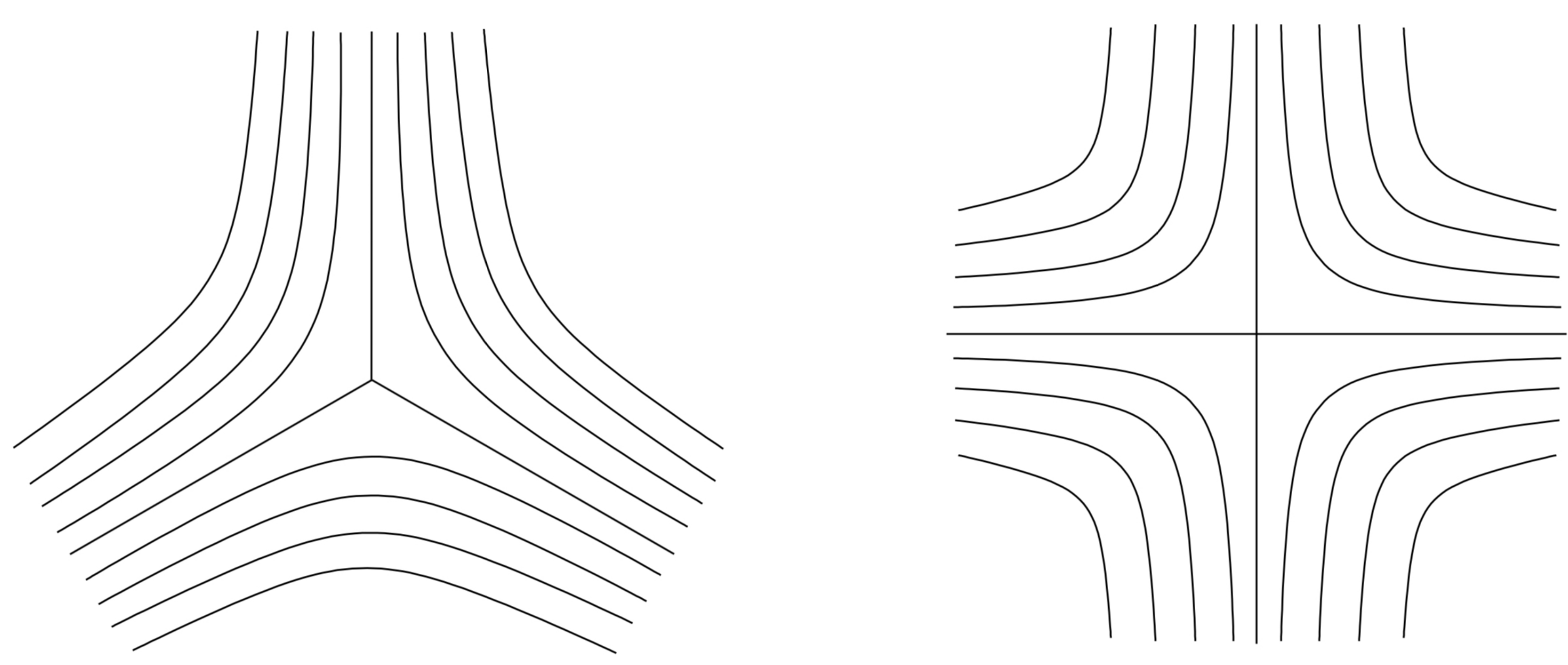}
\caption{Two possible singularities for a foliation. Image reproduced with permission from \cite{FarbMargalit:Primer}.}
\label{F:Sing}
\end{figure}

Because $z\mapsto \pm z+C$ preserves the vertical and horizontal foliations of $\bC$, we find that each quadratic differential defines vertical and horizontal foliations called $h(q)$ and  $v(q)$ on the surface. These foliations are singular at the zeros (cone points) of the quadratic differential, as in Figure \ref{F:Sing}. They also come equipped with some extra structure called a transverse measure, which assigns to each arc on the surface a non-negative real number measuring the extent to which the arc crosses the foliation.  (So an arc contained in a leaf of the foliation has 0 transverse measure, and the measure of any arc does not change when the arc is pushed along leaves of the foliation.)  A foliation equipped with a transverse measure is called a measured foliation. 

Let $\MF$ be the space of measured foliations on a surface of fixed genus, up to a natural notion of equivalence. 
Let $\cQ \cT_{g,n}$ be the bundle of non-zero quadratic differentials over Teichm\"uller space, and let $\cQ \cM_{g,n}$ be the corresponding bundle over moduli space. A foundational result gives that the map 
$$\cQ \cT_{g,n}\to \MF\times \MF, \quad\quad (X,q) \mapsto (h(q), v(q))$$
is a homeomorphism onto its image. The complement $\Delta$ of the image admits an explicit description. Thus, given a pair of measured foliations $(h,v)$ not in $\Delta$, we can construct a quadratic differential $q(h,v)$ with these horizontal and vertical measured foliations. 

Let $X$ be a hyperbolic surface. A geodesic lamination on $X$ is a closed subset of $X$ that is a union of disjoint simple geodesics. The simplest example is a simple closed geodesic. We let $\ML$ denote the space of measured geodesic laminations on $X$. The measure refers again to a transverse measure, and rather than giving a definition we mention that, in the case of a union of disjoint closed  geodesics, the data of the measure is equivalent to the data of a positive weight for each of the disjoint closed  geodesics. So if $\alpha$ and $\beta$ are disjoint closed geodesics, $\alpha+\beta$ and $3\alpha+7\beta$ give different points of $\ML$. 

Every measured foliation can be ``tightened" to a measured lamination, roughly by replacing each leaf of the foliation with a corresponding geodesic. This gives a homeomorphism $\MF\to \ML$. Since the space of geodesic laminations is homeomorphic to the purely topological object $\MF$, it doesn't matter exactly which hyperbolic metric is used to define $\ML$. 

For more details, see, for example, \cite[Section 2]{Wright:Earthquake}. 

\subsection{Dynamics.} Given a group action on a space, an invariant measure is called ergodic if it is not the average of two non-proportional invariant measures. Thus, ergodic measures are the indecomposable building blocks of all invariant measures. 

Fix a connected component $\cQ$ of a stratum $\cQ(\kappa)$, and let $\cQ_1\subset \cQ$ denote the subset of unit area surfaces. The locus $\cQ_1$ carries a natural Lebesgue class measure called the Masur-Veech measure, which both Masur and Veech proved has finite total mass. 

A foundational result, proven in the 80s, is that the action of $g_t$ on $\cQ_1$ is ergodic \cite{Masur:Interval, Veech:Geodesic}. We recommend \cite[Section 4]{ForniMatheus:Intro} for an expository account of a proof using modern tools.  Ergodicity here is equivalent to the fact that almost every $g_t$-orbit is equidistributed.

A corollary, which was originally due to Masur \cite{Masur:Ergodic} and also follows very easily from the Mautner Lemma (see, for example, \cite[Lemma 3.6]{BekkaMayer:Ergodic}), is that the action of 
$$u_t=\left( \begin{array}{cc} 1&t\\0&1\end{array} \right)$$ 
on $\cQ_1$ is also ergodic. 

Another corollary, which follows from a general result called the Howe-Moore Theorem (see, for example, \cite{BekkaMayer:Ergodic}), is that the action of $g_t$ is not just ergodic but mixing.  This means that, not only do typical  orbit segments $\{g_t(X,q) : 0\leq t \leq T\}$ equidistribute as $T\to\infty$, but if one considers a nice positive measure set $S$, then the sets $g_T(S)$  equidistribute as $T\to\infty$.

\subsection{Hyperbolicity.}\label{SS:Hyperbolicity} Consider a quadratic differential $(X,q)$, presented using polygons in the complex plane. 
Nudge these polygons in such a way that the real part of each edge vector stays the same, but the imaginary part changes slightly, to obtain a new quadratic differential $(X', q')$. The ``difference" between $(X,q)$ and $(X', q')$ is purely in the imaginary direction, which is contracted by the $e^{-t}$ in
$$g_t = \left( \begin{array}{cc} e^t & 0\\ 0 & e^{-t} \end{array}\right).$$
For this reason, one might hope the distance between $g_t(X,q)$ and $g_t(X',q')$ decays like $e^{-t}$ as $t\to\infty$. 
This naive hope is dashed by the issue of cut and paste, but nonetheless Forni showed that typically the distance decays like $O(e^{-ct})$ for some $0<c<1$ \cite{Forni:Deviations}. See the survey \cite{ForniMatheus:Intro} for more details. 

This contraction effected by the flow $g_t$ is a characteristic feature of geodesic flows on negatively curved manifolds, and adds to dynamical similarities previously established by Veech and others between these two situations \cite{Veech:Geodesic}. 

%
%

\section{Earthquake flow}\label{S:Earthquake}

In this section we describe the remarkable bridge Mirzakhani built in \cite{Mirzakhani:Earthquake} between hyperbolic and flat geometry. The author has already written a survey devoted solely to this topic, which the reader can consult for more details \cite{Wright:Earthquake}. 

\subsection{The definition of earthquake flow.} For each $\lambda\in \ML$, there is a map 
$$E_\lambda: \cT_g\to \cT_g$$
called the earthquake in $\lambda$. Earthquake flow is defined as the family of maps, defined for $t\in \bR$,
$$E_t : \ML\times \cT_g \to \ML\times \cT_g, \quad\quad (\lambda, X)\mapsto (\lambda, E_{t\lambda} X).$$

Earthquake flow is most easily defined for multi-curves $\lambda=\sum_{i=1}^k c_i \gamma_i$. In this case, we can take a pants decomposition that contains all the curves $\gamma_i$, and consider the associated Fenchel-Nielsen coordinates, which consist of a length coordinate and a twist coordinate for every curve in the pants decomposition. Then $E_\lambda(X)$ is defined as the result of adding $c_i$ to the twist coordinate corresponding to $\gamma_i$, and leaving the other coordinates unchanged.  
In other words, $$E_\lambda(X)=\Tw_{\gamma_k}^{c_k} \circ\cdots\circ \Tw_{\gamma_1}^{c_1}(X)$$ can be obtained from $X$ by cutting along each $\gamma_i$ and re-gluing it with a twist of $c_i$. Here we are using the notation of Section \ref{SS:twist}, except that we are omitting the marking. 

Multi-curves are dense in $\ML$, and the earthquake in a general lamination is defined by continuity: If $\lambda\in \ML$ is a limit of multi-curves $\lambda_n$, then we define $E_\lambda(X)$ as $\lim_{n\to \infty} E_{\lambda_n}(X)$. It isn't obvious, but it turns out that this is well-defined, in that the limit is the same even if one uses a different sequence $\lambda_n'$ converging to $\lambda$. 

Earthquake flow descends to a flow on the bundle $\cP \cM_g$ of measured laminations over moduli space. Its study is motivated by its naturality and its applications. 
\begin{itemize}
\item In addition to the geometric definition above, earthquake flow arises as a Hamiltonian flow. 
\item Theorems about earthquake flow, like that any two points of Teichm\"uller space can be joined via an earthquake path $E_{t\lambda}(X), t\in \bR$, or that each length function $\ell_\gamma$ is convex along each earthquake path, are broadly useful in Teichm\"uller theory \cite{Kerckhoff:Nielsen}. 
\item  Sections \ref{S:Ergodic} and \ref{S:Orbits} rely on the results that we turn to now. 
\end{itemize}

\subsection{A measurable conjugacy.} Mirzakhani relates earthquake flow to  part of the $\GL$ action on the space of quadratic differentials. 

\begin{thm}\label{T:main}
There is a measurable conjugacy $F$ between the earthquake flow $E_t$ on $\ML \times \cT_g$ and 
the  action of 
$$u_t=\left(\begin{array}{cc} 1& t\\ 0 & 1\end{array}\right)$$
on $\cQ \cT_g$.
\end{thm}

That $F$ is a measurable conjugacy means that $u_t \circ F = F \circ E_t$ and that $F$ is a measurable bijection between full measure subsets of $\ML \times \cT_g$ and $\cQ \cT_g$.  Theorem \ref{T:main} is surprising in light of known differences between horocycle and earthquake flow paths \cite{Fu:Cusp}, \cite[Proposition 8.1]{Minsky-Weiss:Nondivergence}.

The pre-image of the locus $\cQ_1 \cT_{g}$ of unit area quadratic differentials is the bundle $\cP_1 \cT_g$ of pairs $(\lambda, X)$ where $\ell_\lambda(X)=1$. Again, we do not define the hyperbolic length function $\lambda \mapsto\ell_\lambda(X)$ directly, except to say that it is uniquely specified by  continuity and being equal to usual hyperbolic length in the case where $\lambda$ is a multi-curve.  

$F$ is equivariant with respect to the action of the mapping class group. The quotient $\cP_1 \cM_g= \cP_1 \cT_g / \MCG$ has a natural invariant measure, defined so that the measure of a set $S\subset \cP_1 \cM_g$ is given by 
$$\int_{\cM_g} \mu_{Th}(\{c \lambda: (X,\lambda)\in S, 0<c<1\}) \dVol_{WP}(X).$$
Mirzakhani used a result of Bonahon and S\"ozen \cite{BonahonSozen:Symplectic} to show that $F$ is measure preserving with respect to the Masur-Veech measure on $\cQ_1 \cT_g$ and the lift of this measure to $\cP_1 \cT_g$.  Since the action of $u_t$ on $\cQ_1 \cM_g=\cQ_1 \cT_g/\MCG$ is ergodic, Mirzakhani obtained the following. 

\begin{cor}
The earthquake flow on $\cP_1 \cM_g$ is ergodic. 
\end{cor} 

Since $F$ is measure preserving, Mirzakhani also concluded that the total volume of $\cQ_1 \cM_g$ is equal to 
$$ \int_{\cM_g} B(X) \dVol_{WP} .$$
Note that the quantity    
$$B(X)=\mu_{Th}(\{\lambda \in \ML: \ell_\lambda(X)\leq 1\}).$$
 also appeared in Theorem \ref{T:Ann}. A single formula intertwines the Thurston measure, the Weil-Petersson measure, and the Masur-Veech measure, and is moreover related to earthquakes and counting simple closed curves! 

\subsection{Horocyclic foliations.}\label{SS:horo} 
\begin{figure}[ht!]
\includegraphics[width=.5\linewidth]{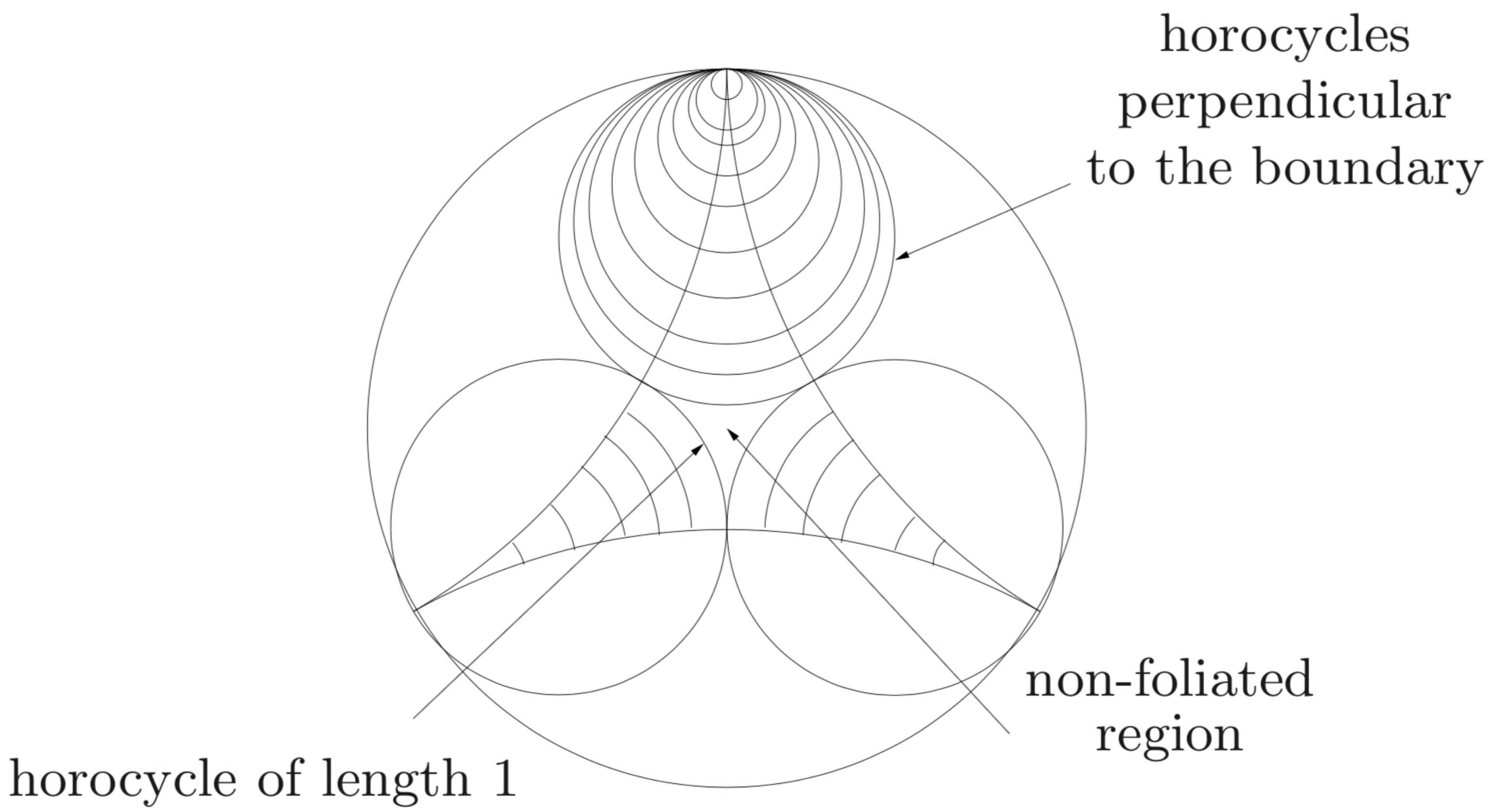}
\caption{Foliating an ideal triangle by horocycles. Shown in the disc model of $\bH^2$. Figure reproduced with permission from \cite{PapadopoulosTheret:Shift}.}
\label{F:foliating}
\end{figure}

A geodesic lamination $\lambda$ is called maximal if $X$ does not contain any bi-infinite geodesics disjoint from $\lambda$. For any such $\lambda$, Thurston defined a map $F_\lambda: \cT_g\to \MF$, which he proved is a  homeomorphism onto its image \cite{Thurston:Stretch}.
Here $F_\lambda(X)$ denotes a specific foliation transverse to $\lambda$ called the horocyclic foliation, defined as follows. The maximal assumption guarantees that $X-\lambda$ is a finite union of ideal triangles. Each can be foliated by horocycles, as in Figure \ref{F:foliating}, leaving a non-foliated region that can be collapsed to a singularity of the foliation. One can modify this foliation on $X-\lambda$  to obtain a foliation of $X$, and one can endow it with a natural transverse measure.  

 Each $\lambda$ can also be viewed as a measured foliation, and for any $X$ there is a quadratic differential 
$$F(\lambda, X) = q(\lambda, F_\lambda(X))$$
whose horizontal foliation is $\lambda$ and whose vertical lamination is $F_\lambda(X)$.  

This map $F$, defined on the full measure set where $\lambda$ is maximal, gives Theorem \ref{T:main}. Mirzakhani shows it is a conjugacy using that both $E_t$ and $u_t$ are Hamiltonian flows, but this can also be seen without reference to the symplectic structure. 

\subsection{Open problems.} It is known that Mirzakhani's measurable conjugacy does not extend to a continuous conjugacy. Remarkably, it isn't known whether some other continuous conjugacy exists. 

\begin{prob}
Is there a continuous conjugacy between earthquake flow and horocycle flow? 
\end{prob}

The author conjectures that there is not. One approach is to try to give a negative answer to the following. 

\begin{prob}
Is the earthquake flow part of a continuous $SL(2,\bR)$ action? 
\end{prob}

We believe that there are variants of Mirzakhani's conjugacy that remain to be explored.

\begin{prob}
Build measurable conjugacies between each stratum of quadratic differentials and certain natural subsets of $\ML\times \cM_g$, as suggested in \cite[Remark 5.6]{Wright:Earthquake}. These conjugacies may extend to loci of quadratic differentials with horizontal saddle connections between the zeros. 
\end{prob}

\begin{prob}
Can some version of the horocyclic foliation $F_\lambda$ be defined for arbitrary non-maximal $\lambda$? 
\end{prob}

Yi Huang suggested that one might try to define $F_\lambda$ so that its leaves contain the level sets of the nearest point projection to $\lambda$. (This projection isn't defined on the geodesic graph of points with more than one closest point on $\lambda$. See \cite[Section 3.2]{Do:Thesis} for some relevant results.)

\section{Horocyclic measures}\label{S:Ergodic}

In this section we give an overview of the papers \cite{LindenstraussMirzakhani:PML} and \cite{Mirzakhani:Random}, which give somewhat related results concerning Teichm\"uller unipotent flow $u_t$ and earthquake flow $E_t$ respectively. 

\subsection{Warm up.} We begin by explaining the theorem of Masur that was used in Section \ref{S:CountSimple}, in order to illustrate how the action of the mapping class group on $\MF$ can be related to the action of unipotent flow on $\cQ_1\cM_{g,n}$. 

\begin{thm}
The action of the mapping class group on $\MF$ is ergodic with respect to the Thurston measure $\mu_{Th}$. 
\end{thm}

This means that any invariant measure $\nu$ that is absolutely continuous with respect to $\mu_{Th}$ is a multiple of $\mu_{Th}$.

Consider any mapping class group invariant measure $\nu$ on $\MF$. Define a measure $\tilde{\nu}$ on the bundle of quadratic differentials $\cQ_1 \cT_{g,n}$, implicitly using the isomorphism between $\cQ_1 \cT_{g,n}$ and a subset of  $\MF\times \MF$, by 
$$\tilde{\nu}(A) = \nu\times \mu_{Th}(\{t q: q\in A, 0<t<1\}),$$
for any $A\subset \cQ_1 \cT_{g,n}$. This $\tilde{\nu}$ has two important properties:  it is $\MCG$-invariant, and it is $u_t$-invariant. The $u_t$-invariance is very important and not hard to prove, but won't be obvious to non-experts. 

The Masur-Veech measure on $\cQ_1 \cT_{g,n}$ is $\tilde{\mu}_{Th}$. It is the pull-back of the Masur-Veech measure on $\cQ_1 \cM_{g,n}$. Recall that the $u_t$ action on $\cQ_1 \cM_{g,n}$ is ergodic. 

If $\nu$ is absolutely continuous with respect to $\mu_{Th}$, then $\tilde{\nu}$ is absolutely continuous with respect to $\tilde{\mu}_{Th}$. One can show that the ergodicity of the $u_t$ action on $\cQ_1 \cM_{g,n}$ implies that $\tilde{\nu}=c\tilde{\mu}_{Th}$ for some $c>0$, and that  that implies $\nu=c\mu_{Th}$. 

\subsection{Ergodic theory on $\MF$.} Mirzakhani and Lindenstrauss classified mapping class group invariant locally finite ergodic measures $\mu$ on $\MF$. To do so, they used that $\tilde{\mu}$ is not just $u_t$-invariant, but it is also horospherical, roughly meaning that it can be studied using the mixing of geodesic flow $g_t$. This connection to mixing is complicated here due to the non-compactness of $\cQ_1 \cM_g$, but this difficulty can be overcome by extending the quantitative non-divergence results for the $u_t$ action proven in \cite{Minsky-Weiss:Nondivergence}. 

For any subsurface $R$, there is a natural inclusion $$\cI_R: \MF(R)\to \MF$$ from the space $\MF(R)$ of measured foliations  on the subsurface to the space of measured foliations on the whole surface.  If $R$ is bounded by closed curves $\gamma_1, \ldots, \gamma_k$, then for any $c=(c_1, \ldots, c_k)\in \bR_+^k$ Mirzakhani and Lindenstrauss consider the map $\cI_R^c: \MF(R)\to \MF$ defined by 
$$\cI_R^c(\alpha)=\cI_R(\alpha)+\sum c_i\gamma_i.$$
If $\mu_{Th}^R$ is the Thurston measure on $\MF(R)$, Mirzakhani and Lindenstrauss observe that $(\cI_R^c)_*(\mu_{Th}^R)$ is a locally finite ergodic measure on $\MF$ invariant under the mapping class group of $S$. Summing over cosets of the mapping class group of $S$ in the mapping class group of the full surface gives a measure 
$$\mu_{Th}^{R,c} = \sum_{[g]} g_* (\cI_R^c(\alpha))$$
that is still locally finite and ergodic, but is now invariant under the whole mapping class group. 

Mirzakhani and Lindenstrauss prove that every $\MCG$-invariant locally finite ergodic measure on $\MF$ is a multiple of $\mu_{Th}$ or $\mu_{Th}^{R,c}$ for some $c$ and $R$. The same result was obtained in \cite{Hamenstadt:Invariant}. 

\subsection{Horospherical measures and earthquake flow.} We now discuss the results of \cite{Mirzakhani:Random}. There is a natural map 
$$\cM_{g-1,2}(L,L)\times \bR/(L\bZ)\to \cM_g$$
obtained by gluing together the two boundary components with a twist given by an element of $\bR/(L\bZ)$. Mirzakhani considers the push-forward of Weil-Petersson measure on $\cM_{g-1,2}(L,L)$ times Lebesgue measure on $\bR/(L\bZ)$, normalized to be a probability measure. She insightfully calls these measures horospherical, but any direct analogy to horospherical measures in other situations remains to be clarified, and her study does not use mixing or  any analogue of the flow $g_t$. Thus the terminology should not be taken literally. 

She shows\footnote{The author and others have not  succeeded in understanding the proof of \cite[Theorem 5.5(b)]{Mirzakhani:Random}, however 
Francisco Arana Herrera has given an alternate proof of \cite[Theorem 5.5(b)]{Mirzakhani:Random}. This is a key step in the proof that the limit measures are absolutely continuous.} 
 that these measures converge  as $L\to\infty$ to a probability measure that is absolutely continuous with respect to the Weil-Petersson measure and has density $B(X)/b_g$, a function that has already appeared in Sections \ref{S:CountSimple} and \ref{S:Earthquake}. This measure is the push-forward via $(X,\lambda)\mapsto X$ of the natural measure on $\cP_1 \cM_g$ that is invariant under the earthquake flow. 

The  map $\cM_{g-1,2}(L,L)\times \bR/(L\bZ)\to \cM_g$ can be lifted to the set of pairs $(X,\gamma)$, where $\gamma$ is the curve obtained from gluing the two boundary components. We can view $\gamma/\ell_\gamma(X)$ as a unit length measured lamination on $X$, and in this way lift the map to the bundle  $\cP_1 \cM_g$ of unit length measured laminations. Mirzakhani uses counting results to show that, in this space, any limit  as $L\to\infty$ of the push-forward measures is absolutely continuous with respect to the natural earthquake flow invariant measure. Any limit must also be invariant under earthquake flow, so she is then able to use ergodicity of earthquake flow to identify   such limit measures up to scale. She then uses non-divergence results from \cite{Minsky-Weiss:Nondivergence}  to show that every limit as $L\to\infty$ is a probability measure, that is, that there is no loss of mass. This concludes the proof that the pushfoward measures converge as $L\to\infty$ to the natural invariant measure on $\cP_1 \cM_g$. 

Other natural methods of building random surfaces in $\cM_g$ are expected to also equidistribute to $B(X)$ times the Weil-Petersson measure. For example, Bainbridge and Rafi informed the author that Mirzakhani advertised the following question. 

\begin{prob}\label{P:pants}
Fix a pants decomposition of a genus $g$ surface. For each $L$, there is a torus in $\cM_g$ obtained by taking each cuff to be length $L$ and looking at all possible twists. Show that these tori equidistribute towards $B(X)$ times the Weil-Petersson measure as $L\to \infty$.
\end{prob}

However, random surfaces with respect to this measure are not yet well-understood. 

\begin{prob}
Study random surfaces in $\cM_g$ sampled with respect to $B(X)$ times the Weil-Petersson measure. For example, what are the asymptotics for the probability that a random surface has a geodesic of length at most $\e$, as $g\to\infty$? 
\end{prob}

\section{Counting with respect to the Teichm\"uller metric}\label{S:ABEM}

In his thesis \cite{Margulis:Anosov}, Margulis considered a compact negatively curved manifold $M=\tilde{M}/\Gamma$ with universal cover $\tilde{M}$ and fundamental group $\pi_1(M)=\Gamma$.  For any $X,Y\in \tilde{M}$, he proved the asymptotic $$|\Gamma\cdot Y \cap B_R(X)|\sim c_{X,Y} e^{hR}$$
for the number of points in  the orbit $\Gamma\cdot Y$ in the ball $B_R(X)$, where $h>0$ is the topological entropy of the geodesic flow, and $c_{X,Y}>0$ is a constant depending on $X$ and $Y$. This is called a ``lattice point count", since one thinks of the orbit $\Gamma\cdot y$ analogously to an orbit of $\bZ^k$ acting on $\bR^k$. He also showed that the number of primitive oriented closed geodesics of length less than $R$ is asymptotic to $\frac{e^{hR}}{hR}.$ 

Instead of such $M$, we wish to consider  $\cM_{g}$ equipped with the Teichm\"uller metric. In comparing $\cM_g$ to a compact negatively curved manifold, one should note the following points. 
\begin{itemize}
\item The Teichm\"uller metric has some features of negative curvature, and is even negatively curved ``on average" \cite{DowdallDuchinMasur:Statistical}. 
\item It is however not actually negatively curved in any sense. For example, there are pairs of infinite geodesic rays in $\cT_g$ leaving the same point that stay bounded distance apart \cite{Masur:Geodesics}. 
\item  There exist closed geodesics outside of any given compact subset of $\cM_{g}$.
\end{itemize} 
Despite these major differences, Mirzakhani and her coauthors obtained  results for $\cM_g$ analogous to those Margulis  obtained for compact negatively curved manifolds \cite{AthreyaBufetovEskinMirzakhani:LatticePoint, EskinMirzakhani:Counting, EskinMirzakhaniRafi:Strata}. 

\subsection{Lattice point counting.} Consider the action of $\Gamma=\MCG$ on $\cT_{g}$, and consider balls $B_R(Y)$ of radius $R$ centered at $Y\in \cT_{g}$. Athreya, Bufetov, Eskin, and Mirzakhani proved the  following. 

\begin{thm}\label{T:ABEM} 
As $R\to\infty$,  
$$|\Gamma\cdot Y \cap B_R(X)|\sim C_{X,Y} e^{hR},$$
where $h=6g-6$ is the entropy of the Teichm\"uller geodesic flow and $C_{X,Y}$ is an explicit constant.
\end{thm}

Dumas and Mirzakhani later discovered that the constant $C_{X,Y}$ does not depend on the choice of $X,Y\in \cT_{g}$ \cite[Theorem 5.10]{Dumas:Skinning}. It would be interesting to understand the asymptotics of this constant as $g\to\infty$. 

In comparison to the problem of counting points in $\bZ^k$ in a ball, the difficulty of this problem is that a definite proportion of the points of $\Gamma\cdot Y \cap B_R(X)$ lie near the boundary sphere of $B_R(X)$; this is witnessed by the exponential growth of the count. 
A solution is to use mixing of the $g_t$ action to prove that large spheres are equidistributed in $\cM_g$, and to apply this equidistribution to understand how many points of $\Gamma\cdot Y$ each concentric sphere of $B_R(X)$ is close to. For an introduction to this technique in the simplest situation, we highly recommend \cite[Section 2]{EskinMcMullen:Mixing}. In this situation, a key technical issue is the smaller dimensional strata of the bundle of quadratic differentials, close to which the flow $g_t$ is less understood. 

Athreya, Bufetov, Eskin, and Mirzakhani also give asymptotics for the volume of $B_R(X)$.

\subsection{Counting closed orbits.}  Let $N(R)$ denote the number of primitive oriented closed geodesics on $\cM_{g}$ of length at most $R$. 

Following Thurston, every element of the mapping class group can be classified as periodic, reducible, or pseudo-Anosov \cite[Chapter 13]{FarbMargalit:Primer}. Briefly, an element is pseudo-Anosov if none of its powers fix a multi-curve; these elements are thought of as ``generic" in the mapping class group. 

Conjugacy classes of pseudo-Anosovs are in bijection with oriented closed geodesics on $\cM_g$, and the translation length of the pseudo-Anosov acting on $\cT_g$ is the length of the closed geodesic. Thus, $N(R)$ can also be interpreted as the number of conjugacy classes of primitive pseudo-Anosovs of translation length at most $R$. 

Eskin and Mirzakhani show the following. 

\begin{thm}\label{T:NR}
As $R\to\infty$, 
$$N(R)\sim \frac{e^{hR}}{hR}.$$
\end{thm}

Later, Eskin, Mirzakhani and Rafi proved a version of this for the $g_t$ action on individual strata of quadratic differentials. 

The basic idea of the proof is to use mixing of the $g_t$ action to count the number of orbit segments of $g_t$ that come back close to where they started, and to prove a Closing Lemma to show that typically each such orbit segment lies close to a closed orbit.  Following a standard argument, the Closing Lemma is deduced from two applications of the contraction mapping principle and the hyperbolicity results discussed in Section \ref{SS:Hyperbolicity}. 

%

This strategy, however, only gives information on the closed geodesics that intersect some fixed compact set in $\cM_g$. And $\cM_g$ contains a great many closed geodesics outside of any compact set \cite{Hamenstadt:Thin}. 

The key to counting closed geodesics in this context is thus to show that there are not too many that live entirely outside a compact set. More precisely, if $K$ is a large compact set, the authors show that the number of closed geodesics of length $R$ that stay outside of $K$ has a smaller exponential order of growth than $N(R)$. 

To prove this, the authors consider a discretization of a geodesic into a sequence of bounded length segments between a discrete ``net" of points in $\cT_g$.  They then show that most such sequence do not avoid a compact set. In a precise sense, a random sequence is biased towards returning to compact sets, and so it is very unlikely that a long random sequence entirely avoids a fixed large compact set.  

This final step is more difficult in the case of strata, and in fact the problem is open in a more general context discussed in Section \ref{S:Measures}. 

\begin{prob}
Prove a version of Theorem \ref{T:NR} that counts closed $g_t$ orbits in a $SL(2,\bR)$--orbit closure in a stratum. 
\end{prob}

\section{From orbits of curves to orbits in Teichm\"uller space}\label{S:Orbits}

 In Mirzakhani's proof of Theorem \ref{T:Ann}, discussed in Section \ref{S:CountSimple}, it was absolutely crucial that $\gamma$ be a simple multi-curve, so that it defines a point of $\MF$, and so that the integration formula can be applied. In one of her most recent papers, which is supplemented by forthcoming work of Rafi and Souto, Mirzakhani gave a totally different approach to this problem that allowed for $\gamma$ to have self-intersections \cite{Mirzakhani:Orbits}. 

\begin{thm}\label{T:LastCount}
Let $\gamma$ be any homotopy class of closed curve.\footnote{Most of \cite{Mirzakhani:Orbits} is devoted to the case of filling curves. Mirzakhani claims that the non-filling case is similar, referring in the last paragraph of the proof of Theorem 1.1 to a remark in Section 6.3. That remark seems incorrect. 
In forthcoming work, Rafi and Souto explain how to deduce the non-filling case from the filling case \cite{RafiSouto:Revisited}.}
Then
$$\lim_{L\to \infty} \frac{s_X(L, \gamma)}{L^{6g-6+2n}} = \frac{c(\gamma) \cdot B(X)}{b_{g,n}},$$
where $B(X)=\mu_{Th}(B_X)$ and $b_{g,n}=\int B(X) \dVol_{WP}$ are as in Sections \ref{S:CountSimple} and \ref{S:Earthquake} and $c(\gamma)\in \bQ^+$. 
\end{thm}

For any $k$, the set of curves with exactly $k$ self-intersections is a finite union of mapping class group orbits, so this result implies asymptotics for the set of curves with $k$ self-intersections. In the case of a once punctured torus, this result can be used to count integer solutions to the Markoff equation $p^2+q^2+r^2=3pqr$ \cite[Section 3]{Mirzakhani:Orbits}. 

In forthcoming work, Erlandsson and Souto will give an independent proof of Theorem \ref{T:LastCount}.

\subsection{Other notions of length.} In the same way that each simple closed curve defines a point in $\MF$, each closed curve defines a point in the space $\cC$ of geodesic currents. You can think of this space as a completion of the space of linear combinations of closed curves, in the same way that $\MF$ is a completion of the space of linear combinations of simple multi-curves. 

One can define the length on $X$ of any current; this is uniquely specified by the fact that the length function $\cC\to \bR_+$ is continuous. More generally, one could consider any reasonably nice function $f:\cC\to\bR_+$ as defining a type of length on currents. There are numerous important examples of such functions $f$, including length on a fixed Riemannian metric with negative but non-constant curvature, and notions of length coming from Higher Teichm\"uller Theory. 

A combination of Mirzakhani's result and work of Erlandsson and Souto  give that Theorem \ref{T:LastCount} is true even if we replace $s_X(L,\gamma)$ with 
$$s_f(L,\gamma) = \left|\{ \alpha \in \MCG\cdot \gamma : f(\alpha)\leq L   \}  \right|,$$
for any homogeneous and continuous $f:\cC\to\bR_+$ \cite{ErlandssonSouto:Counting}. Roughly speaking, Erlandsson and Souto show that the problem of counting $s_f(L,\gamma)$ does not depend on which $f$ is used, and Theorem \ref{T:LastCount} solves this problem when $f$ is the usual  length on some hyperbolic surface $X$.   

 Rafi and Souto additionally showed that one can count mapping class group orbits of more general currents, and also identified Mirzakhani's constant $c(\gamma)$ \cite{RafiSouto:Currents}. See  \cite{ErlandssonUyanik:Length} for a survey of related developments.

\subsection{The idea of the proof.} Take $\gamma$ a filling curve. By definition, this means that $\gamma$ intersects every simple closed curve. 

Mirzakhani transforms the problem to a counting problem on Teichm\"uller space. First, she observes that $\Stab(\gamma)$ is finite and
$$  s_X(L,\gamma)=  
\frac{\left|\{ g\in \MCG : \ell_X(g\cdot \gamma)\leq L   \}  \right|}{\left| \Stab(\gamma)\right|}.$$
But $\ell_X(g\cdot\gamma) = \ell_{g^{-1}\cdot X}(\gamma)$, and $\Stab(X)$ is finite for any $X$, so we get 
\begin{eqnarray*}
&&\left|\{ g\in \MCG : \ell_X(g\cdot \gamma)\leq L   \}  \right| 
\\&=& 
\left|\{ g\in \MCG : \ell_{g\cdot X}( \gamma)\leq L   \}  \right|
\\&=&
\left| \Stab(X)\right|\cdot\left|\{ Y\in \MCG\cdot X : \ell_{Y}( \gamma)\leq L   \}  \right|.
\end{eqnarray*}
In summary,  the problem is reduced to counting the number of elements in the mapping class group orbit of $X$ that lie in the compact set 
$$B_\gamma(L) = \{Y\in \cT_{g,n}: \ell_\gamma(Y)\leq L\}. $$ 
This problem is  similar to the lattice point count of the previous section, except that the shape of the set $B_\gamma(L)$ is a priori not as well understood. 

Mirzakhani shows that in Fenchel-Nielsen coordinates, $B_\gamma(L)$ isn't too strangely shaped: it is asymptotically polyhedral. In the previous lattice point counting problem, Mirzakhani and her coauthors used equidistribution of certain spheres in $\cM_g$. Similarly, here Mirzakhani uses that the boundary of $B_\gamma(L)$  equidistributes in moduli space as $L\to\infty$. She does this by showing that the boundary can be divided into pieces that can be approximated by pieces of  ``horospheres" similar to those studied in Section \ref{SS:horo}, and using equidistribution of those ``horospheres".

\section{$SL(2,\bR)$--invariant measures and orbit closures}\label{S:Measures}
McMullen classified $SL(2,\bR)$-invariant measures and orbit closures of Abelian differentials in genus 2 \cite{McMullen:Dynamics}; see also \cite{Calta} for related work. Mirzakhani was deeply interested in generalizing these results  to higher genus. In this section, we discuss the results from \cite{EskinMirzakhani:InvariantMeasures, EskinMirzakhaniMohammadi:Isolation} on this topic. Lecture notes of Eskin, a short survey by the author, and an in-depth survey by Quint all give more details \cite{Eskin:Notes, Wright:Billiards, Quint:Rigidite}. 

First, Eskin and Mirzakhani understood 
$$P= \left\{\left(\begin{array}{cc} * & *\\ 0& * \end{array}\right)\right\}\subset SL(2,\bR)$$
invariant measures on moduli spaces of translation surfaces.

\begin{thm}\label{T:EM}
Every $P$-invariant ergodic probability measure on the locus of unit area surfaces in a stratum is the natural Lebesgue measure on the unit area locus of a (properly immersed) linear sub-orbifold of the stratum. In particular, it is $SL(2,\bR)$-invariant. 
\end{thm}

A sub-orbifold is linear if it is locally described by linear equations on the edges of polygons defining the surfaces, and if these equations have zero constant term and real coefficients.  

Then, Eskin, Mirzakhani and Mohammadi used the previous result to prove the following.

\begin{thm}\label{T:EMM}
The $P$-orbit closure of any point is equal to the $SL(2,\bR)$-orbit closure of that point, and is the unit area locus of a linear sub-orbifold of the stratum. 
\end{thm}

Before giving the definitions, we give some of the reasons why so much interest has focused on these questions, and why these theorems have been called a ``Magic Wand" for applications \cite{Zorich:Magique, Zorich:Magic}. 

\begin{itemize}
\item Since each $SL(2,\bR)$ orbit projects to a complex geodesic in $\cM_{g,n}$, orbit closures are related to the geometry of $\cM_{g,n}$. 
\item Quadratic differentials arise naturally in the study of billiards in rational polygons. Questions such as the number of generalized diagonals in a polygon can be productively studied using, and in fact depend on, the $SL(2,\bR)$ orbit closure. 
\item There are  applications to other simplified models in physics, such as the Ehrenfest Wind Tree Model and periodic arrays of Eaton lenses \cite{FraczekShiUlcigrai:Genericity}. 
\item There are  applications to other low complexity dynamical systems such as interval exchange transformations and flows on surfaces, as well as the illumination problem \cite{LelievreMonteilWeiss:Everything}. (If the edges of a non-convex polygon were reflective, what subset of the polygon would be illuminated by a light source in the polygon?)
\item The orbit closures  turn out to be beautiful and independently interesting algebraic varieties. 
\item There are deep connections to homogeneous space dynamics, which have enriched both fields. 
\end{itemize}

 Each point of the $SL(2,\bR)$-orbit of a surface $(X,q)$ represents a different perspective on $(X,q)$, which in a sense arises from a linear change of coordinates. Understanding the orbit of a surface represents understanding all possible perspectives on the surface simultaneously.

\subsection{Definitions.} We begin with a warm up observation. Consider the diagonal action of $GL(2,\bR)$ on $\bC^m \simeq (\bR^2)^m$. A subspace of $\bC^m$ is invariant under this action if and only if it can be cut out by linear equations with real coefficients. This explains the definition of linear submanifold. It is almost as easy to see that a linear submanifold of a stratum must be invariant by $GL(2,\bR)$, since this can be checked locally. 

Each linear submanifold $\cM$ inherits from the stratum an atlas of charts to $\bC^m$ with transition functions in $GL(m,\bR)$. Suppose the transition functions can be chosen to have determinant $\pm1$. For Theorem \ref{T:EM} it suffices to deal with this case, and actually it turns out this is always possible. Then there is a natural Lebesgue $\mu$ measure on  $\cM$, well-defined up to scale.  From this, one can define a measure $\mu_1$ on the unit area locus by 
$$\mu_1(A) = \mu(\{t (X,\omega) : t\in (0,1), (X,\omega)\in A    \}).$$
It is exactly these measures that we refer to in Theorem \ref{T:EM}. 

\subsection{Previous results in measure rigidity.} Theorems \ref{T:EM} and \ref{T:EMM} were inspired by and built upon related results in homogeneous space dynamics, where one considers group actions on homogeneous spaces rather than moduli spaces of surfaces. 

The philosophy behind this line of work is that, even if one wants primarily to understand orbit closures, it is more convenient to first work with invariant measures. For a short introduction to this philosophy, see Venkatesh's survey \cite{Venkatesh:EKL}.   

Any student seriously interested in learning the proofs of Theorems \ref{T:EM} and \ref{T:EMM} might be well advised to start in the homogeneous setting, for example by reading about Ratner's Theorems from expository sources  such as \cite{Einsiedler:Ratner, Eskin:Unipotent, Morris:Ratner}, or by reading the work of Benoist and Quint \cite{BenoistQuint:Mesures, Translation}. Let us briefly mention the work of Benoist and Quint in the simplest example. 

\begin{thm}\label{T:BQ}
Let $\mu$ be a finitely supported measure on $SL(2,\bZ)$. If $\nu$ is any non-atomic $\mu$-stationary probability measure on $\bR^2/\bZ^2$, then $\nu$ is Lebesgue. 
\end{thm}

Here $SL(2,\bZ)$ acts on the torus $\bR^2/\bZ^2$ in the usual way, and $\nu$ is called $\mu$-stationary if it is ``invariant on average" by the $SL(2,\bZ)$ action: 
$$\nu = \sum_{g\in SL(2,\bZ)} \mu(g) g_* \nu.$$

In particular, every invariant measure is stationary, but there are many important situations where one cannot prove the existence of invariant measures but can easily show that stationary measures exist. Stationary measures are closely related to the study of random walks.

Using Theorem \ref{T:BQ}, Benoist and Quint are able to also understand orbit closures, and show that any closed  set invariant under $SL(2,\bZ)$ must be finite or the whole torus $\bR^2/\bZ^2$. 

In addition to the ideas of Ratner, Benoist, and Quint, and new ideas, Eskin and Mirzakhani make use of ideas developed in the context of homogeneous dynamics by mathematicians including Margulis, Tomanov, Einsiedler, Katok, and Lindenstrauss \cite{MargulisTomanov:Invariant, EinsiedlerKatokLindenstrauss:Invariant}.

\subsection{Structure of the proofs.} The techniques used by Eskin, Mirzakhani, and Mohammadi are very abstract, and do not directly involve  the geometry of surfaces. 

Instead, Eskin and Mirzakhani focus on entropy and the rates at which the geodesic flow $g_t$ expands and contracts. These rates are called Lyapunov exponents, and are the subject of work of Forni and others in this context \cite{ForniMatheus:Intro}. They also use an idea of Furstenburg to convert the problem about $P$-invariant measures to a problem about measures that are stationary with respect to certain measures on $SL(2,\bR)$. 

The hardest and longest part of the proof of Theorem \ref{T:EM} works in this framework to prove that a stationary measure is in fact $SL(2,\bR)$-invariant. Once the measure is known to be $SL(2,\bR)$-invariant, Eskin and Mirzakhani use a version of the arguments of Benoist and Quint to conclude. The general framework of the proof is to establish additional invariance properties bit by bit, but this is complicated by the fact that here the additional ``invariance" is often not given by a global action of a group. 

The algebraic properties of $P$, namely that it is amenable, allow Eskin, Mirzakhani and Mohammadi to  ``average" over partial $P$-orbits to construct, from each $P$-orbit, a $P$-invariant measure. This would not be possible with $SL(2,\bR)$, and is the chief reason that Eskin and Mirzakhani worked so hard to understand $P$-invariant measures, instead of the easier $SL(2,\bR)$-invariant measures. 

Eskin, Mirzakhani and Mohammadi show that the support of the measure constructed from each orbit is equal to that orbit closure,  and in this way prove Theorem \ref{T:EMM}.  

\subsection{Subsequent work and open problems.} The techniques of \cite{EskinMirzakhani:InvariantMeasures, EskinMirzakhaniMohammadi:Isolation} have since been applied to homogeneous and smooth dynamics \cite{EskinLindenstrauss:Zariski, EskinLindenstrauss:Random, BrownRodriguezHertz:Rigidity}.

Mirzakhani was deeply interested in the $u_t$ action. 

\begin{prob}
Understand $u_t$-invariant measures and orbit closures. 
\end{prob}

The hope was originally for a version of Ratner's Theorems in this context. However, a recent discovery of Chaika, Smillie and Weiss shows that $u_t$-orbit closures can have non-integral Hausdorff dimension, so the theory for homogeneous spaces cannot hold in full in this context \cite{ChaikaSmillieWeiss}. However, one can still hope for progress for particular classes of $u_t$-invariant measures. 

\begin{prob}
Understand piecewise linear $u_t$-invariant measures and orbit closures. 
\end{prob}


We conclude with a question recommended to us by McMullen, concerning a generalization of flat surfaces arising by tensoring the bundle of quadratic differentials on a Riemann surface with a flat line bundle with holonomy in $\bR_+$.

\begin{prob}
What can one say about the action of $SL_2(\bR)$ on bundles $Q'\cM_g\to \cM_g$ of twisted quadratic differentials?
\end{prob}

Twisted quadratic differentials are examples of dilation surfaces, whose  moduli spaces may have features in common with infinite volume hyperbolic manifolds. See \cite[Section 8]{Ghazouani:Tori} for more details, as well as for some conjectures and open problems. 

\section{Classification of $SL(2,\bR)$--orbit closures}\label{S:OrbitClosures}

It turns out that orbit closures of translation surfaces are  varieties that can be characterized by unlikely algebro-geometric properties \cite{Filip:Splitting}, and there are at most countably many of them.  Mirzakhani  proposed the following in answer to the mystery of why so few examples of orbit closures are known. 

\begin{conj}
Every orbit closure of rank at least 2 is trivial. 
\end{conj}

We define these terms for orbit $GL^{+}(2,\bR)$--orbit closures  of Abelian differentials. One can get corresponding definitions for orbit closures of quadratic differentials by passing to double covers where the quadratic differentials are squares of Abelian differentials. The change from $SL(2,\bR)$ to $GL^{+}(2,\bR)$ is unimportant but simplifies the notation. 

  The tangent space $T_{(X,\omega)}\cM$ to an orbit closure $\cM$ of Abelian differentials is naturally a subspace of the relative cohomology group $H^1(X,\Sigma, \bC)$, where $\Sigma$ is the set of zeros of $\omega$.  There is a natural map $p:H^1(X,\Sigma, \bC)\to H^1(X,\bC)$, and 
$$\rank(\cM)=\frac12 \dim_{\bC} p(T_{(X,\omega)} \cM)$$
for any $(X,\omega)\in \cM$. This is an integer between 1 and the genus of $X$. 

\begin{ex}
The rank of a closed orbit is 1. The rank of a stratum of genus $g$ Abelian differentials is $g$. 
\end{ex}

An alternative definition of rank is provided by Filip's work: There is a number field $\bk(\cM)$ that acts naturally on the space  $H^{1,0}(X,\bC)$ of Abelian differentials  for each $(X,\omega)\in \cM$. (In fact $\bk(\cM)$ is totally real, and the action comes from an action on $\Jac(X)$.) The cohomology class of $\omega$ is a eigenvector for this action, and the rank is the dimension of the eigenspace. 

An orbit closure is called trivial if it arises via a branched covering construction from a stratum of Abelian or quadratic differentials. 

\begin{ex}
Consider a connected component $\cM$ of the locus of $(X,\omega)$ for which there exists a $(X', \omega')\in \cH(2)$ and degree 3 cover $f:X\to X'$, simply branched over one point that isn't a zero of $\omega'$, such that $\omega=f^*(\omega')$. This $\cM$ is trivial, since it arises from $\cH(2)$ from a covering construction. (None of the choices of degree or branching behavior is important.)
\end{ex}

Mirzakhani's conjecture is  known to be true in genus 3 \cite{AulicinoNguyen:Three} and in the hyperelliptic case \cite{Apisa:Hyp}, and it is known that there are at most finitely many counterexamples in each genus \cite{EskinFilipWright:Hull}. It is also of true in genus 2, where McMullen's classification preceded Mirzakhani's conjecture \cite{McMullen:Dynamics}. 

However,  8 counterexamples have been discovered, all of rank 2 \cite{McMullenMukamelWright:Cubic,EskinMcMullenMukamelWright:Billiards}. These arise from algebro-geometric constructions and have surprising implications for billiards in certain quadrilaterals. A year of discussions with Mirzakhani, aimed at proving her conjecture, contributed to the discovery of one of the counterexamples.  

\begin{figure}[ht!]
\includegraphics[width=.5\linewidth]{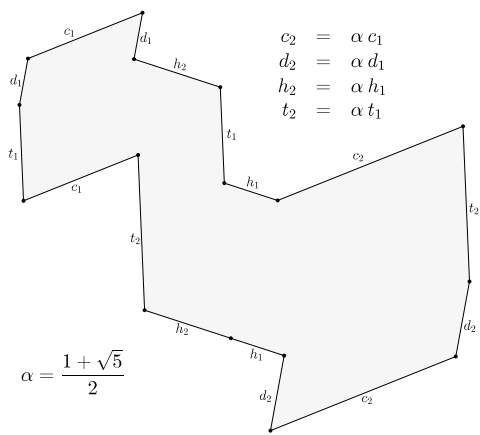}
\caption{A local description of one of the counterexamples.}
\label{F:H6}
\end{figure}
 
Despite the counterexamples, Mirzakhani still wanted to find a strong statement like the conjecture that is true, and to at least partially resolve the mystery of why so few orbit closures are known. To this end, the author and Mirzakhani wrote a paper giving the basic boundary theory of orbit closures, allowing inductive arguments on genus \cite{MirzakhaniWright:Boundary}, see also \cite{ChenWright:WYSIWYG}. This was later used in \cite{MirzakhaniWright:FullRank} to prove the conjecture in the case of maximum rank, which is rank equal to the genus. This also used the author's Cylinder Deformation Theorem from \cite{Wright:Cylinder}, which says that often cylinders in translation surfaces can be deformed while staying in the orbit closure. These deformations can be used to degenerate, facilitating inductive arguments.  

We believe that the most important open problem on orbit closures is now the following. 

\begin{prob}
Is every orbit closure of rank at least 3 trivial? 
\end{prob}

Some limited theoretical and computer efforts have failed to find any non-trivial orbit closures of rank at least 3 thus far. Furthermore, every $(X,\omega)$ with a rank 2 orbit closure has a natural map $X\to\bP^1$, which played a crucial role in the 8 counterexamples. However, in the rank 3 case there is instead a map $X\to\bP^2$, and so far there is no proposal for how this map could be leveraged to establish the special algebro-geometric properties that, by the work of Filip, must be present.

There are myriad other  open problems about orbit closures, but we conclude by mentioning just two that Mirzakhani was interested in. 

\begin{prob}
Give an effective proof that there are only finitely many orbit closures of rank at least 2 in each genus. 
\end{prob}

\begin{prob}
Classify algebraically primitive orbit closures of rank at least 2, where an orbit closure is algebraically primitive if $\rank(\cM)\cdot [\bk(\cM):\bQ]=g$. Start with the case where $\bk(\cM)$ is quadratic. 
\end{prob}

See \cite{Wright:Broad} for the author's more in depth (although now somewhat out of date) survey aimed towards orbit closures, and see \cite{MirzakhaniWright:FullRank} for a few  additional open problems.

\section{Effective counting of simple closed curves}\label{S:Effective}

We now discuss work of Eskin, Mirzakhani, and Mohammadi that gives an effective version of Theorem \ref{T:Ann} with a totally different proof \cite{EskinMirzakhaniMohammadi:Effective}. Recall that if $X$ is a hyperbolic surface and $\gamma$ is a multi-curve, 
$$s_X(L,\gamma) = |\{\alpha \in \MCG \cdot \gamma : \ell_\alpha(X) \leq L\}|.$$

\begin{thm}\label{T:Effective}
There is a constant $\kappa=\kappa(g)>0$ such that for any $X\in \cM_g$ and any rational multi-curve $\gamma$, 
$$s_X(L, \gamma) = \frac{c(\gamma) \cdot B(X)}{b_{g}} L^{6g-6} + O(L^{6g-6-\kappa}).$$
\end{thm}

Here the constants  $c(\gamma), B(X), b_g$ are as in Theorem \ref{T:Ann}, but will not play a major role in our discussion. 

The proof of Theorem \ref{T:Effective} uses the dynamics of Teichm\"uller geodesics flow, and the effective nature of the count comes from the exponential mixing of this flow proven in \cite{AvilaGouezelYoccoz:Exponential}. 

\subsection{A special case.} Rather than discuss the rather technical proof, we illustrate the relationship between the counting problem and geodesic flow in an extremely special and well-known case. Namely, we will outline an effective count of the number of elements of  
$$\bZ^2_{prim} = \{(p,q)\in \bZ^2 : \gcd(p,q)=1\} = SL(2,\bZ) \cdot (1,0)$$ 
of size at most $L$.  

Consider the unit tangent bundle $T^1 \bH^2$ to the hyperbolic plane, viewed in the upper half plane model. $PSL(2,\bR)$ acts on $\bH^2$ by isometries, inducing an action on $T^1 \bH^2$. 

An example of a horosphere in $T^1 \bH^2$ is the set 
$$\left\{ \left( \begin{array}{cc} 1 & t\\ 0 & 1\end{array}\right) (i,i): t\in \bR \right\},$$
where $(i,i)$ denotes the vertical unit tangent vector at the point $i$ in the upper half plane. 
\begin{figure}[ht!]
\includegraphics[width=.9\linewidth]{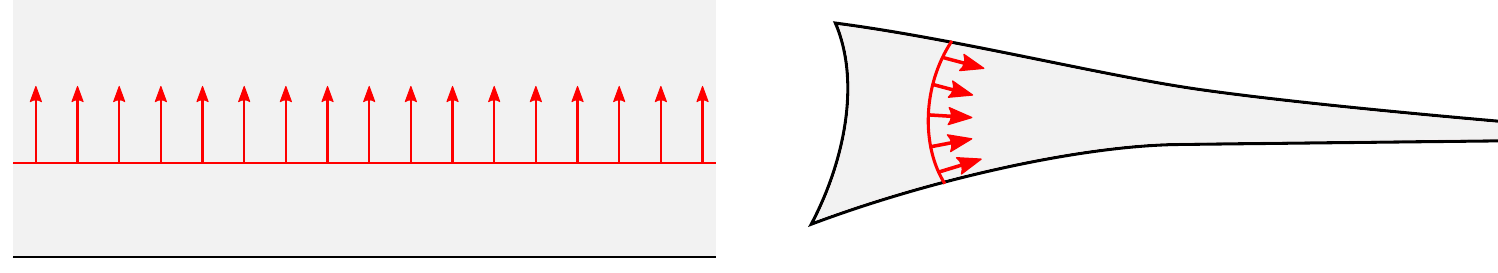}
\caption{A horosphere in $T^1\bH^2$ and its image in $T^1 \bH^2 / SL(2,\bZ)$.}
\label{F:horo}
\end{figure}
Every other horosphere is the image of this horosphere via an isometry, and the stabilizer of this horosphere is 
$$U=\left\{\left( \begin{array}{cc} 1 & t\\ 0 & 1\end{array}\right) \right\}.$$
Thus, we can identify the set of horospheres as $PSL(2,\bR)/U$. 

$PSL(2,\bR)$ also acts on $(\bR^2-\{(0,0)\})/\pm1$, and the stabilizer of $(1,0)$ is $U$. Thus, we can identify $(\bR^2-\{(0,0)\})/\pm1$ with the set of horospheres in $T^1\bH^2$. 

Via this identification, the norm of an element of $\bR^2-\{(0,0)\}$ is a function of the distance of the corresponding horosphere to the point  $i\in \bH^2$ fixed by $SO(2)\subset PSL(2,\bR)$. 

Thus, the original problem concerning $\bZ^2_{prim}$ can be translated to the problem of counting the number of elements of the $SL(2,\bZ)$ orbit of a horosphere that intersect a large ball in $\bH^2$. All these horospheres give rise to a single compact subset of $T^1 \bH^2 / SL(2,\bZ)$ called a closed horosphere. The original problem is thus translated to counting the number of geodesics of length at most $f(L)$ from a fixed point in $\bH^2 / SL(2,\bZ)$ to a fixed closed horosphere, where $f(L)$ is an explicit function. This problem manifestly involves geodesic flow, and ideas of Margulis can be used to solve it effectively. 


\subsection{The hyperbolic length function.} The above outline could be used to count elements of $\bZ^2_{prim}$ of ``size" at most $L$, where the notion of ``size" is given by a  sufficiently nice function. 

Similarly, Eskin-Mirzakhani-Mohammadi must check  that their notion of ``size" is sufficiently nice. They do this with a result from Mirzakhani's thesis \cite[Theorem A.1]{Mirzakhani:Thesis}, which states that the hyperbolic length function 
$$\alpha\mapsto \ell_\alpha(X)$$
is a convex function on natural subsets of $\ML$ called train track charts.   A corollary is that this function is Lipschitz.

We end with the natural next step in this line of investigation. 

\begin{prob}
Generalize Theorem \ref{T:Effective} to non-simple multi-curves, thus giving an effective version of Theorem \ref{T:LastCount}. 
\end{prob}

\section{Random walks on the mapping class group}\label{S:RW}

Consider a reasonably nice measure $\mu$ on the group of isometries of $\bH^2$. One of the archetypal features of negative curvature is that a random walk on $\bH^2$ with steps determined by $\mu$ will converge almost surely to the circle $\partial \bH^2$ at infinity. 

If one fixes a starting point, there is a ``hitting measure" $\nu$ on $\partial \bH^2$, whose defining property is that  the measure $\nu(A)$ of any subset $A\subset \partial \bH^2$ is equal to the probability that a random walk path will converge to a point in $A$. The hitting measure is also sometimes called the harmonic measure, because it satisfies an equation  saying that it is unchanged by averaging over $\mu$. If the support of $\mu$ is discrete, this equation can be written as 
\begin{equation}\label{E:harmonic}
\nu = \sum_\gamma \mu(\gamma) \gamma_*(\nu).
\end{equation} 

Moreover, almost every random walk path tracks a geodesic ray, and the endpoint of the geodesic ray is the point in $\partial \bH^2$ that the random walk path converges to. More precisely, the geometric form of the Oseledets Multiplicative Ergodic Theorem states that the distance between the   the random walk path and the geodesic ray that it tracks is sub-linear in the number of steps of the random walk. 

A natural way to pick a random geodesic ray leaving a point is to simply pick the endpoint of the ray using the Lebesgue measure on the circle. This notion captures what is commonly understood as a random geodesic ray. 

One may ask: How different are random geodesic rays and random walk paths? Can one simulate the effect of picking a random geodesic ray by instead generating a random walk path? This question has extra significance if we want the the measure $\mu$ defining the random walk to live in a discrete subgroup such as $SL(2,\bZ)$. Can one simulate a random geodesic ray on $\bH^2/\SL(2,\bZ)$ by a random walk on its fundamental group $SL(2,\bZ)$? 

At first, this seems like a lot to ask for. Indeed, although the geodesic flow on the unit tangent bundle $T^1 \bH^2/\SL(2,\bZ)$ is very chaotic, it is of course a deterministic process. In contrast, the random walk has  no ``memory" at all, in that each step of the random walk does not depend at all on the previous steps. For this reason, random walks are easier to study than the geodesic flow.

Nonetheless, Furstenberg famously proved that one can simulate the random geodesic ray with a random walk path: There is a measure on $SL(2,\bZ)$ whose hitting measure on $\partial \bH^2$ is Lebesgue. Furstenberg obtained his measure $\mu$ as a discretization of Brownian motion, and used it to prove one of the first rigidity results, for example that groups like $SL(n,\bZ), n>2$  cannot be isomorphic to  subgroups of $SL(2,\bR)$. For a very clear exposition, see \cite[Sections 5, 6]{Furstenberg:Walks}.

In forthcoming work, Eskin, Mirzakhani and Rafi will give an analogue of Furstenberg's result in the context of Teichm\"uller theory \cite{EskinMirzakhaniRafi:Walks}. Thus $\bH^2$ will be replaced with Teichm\"uller space, and the circle $\partial \bH^2$ at infinity will be replaced with the Thurston boundary of Teichm\"uller space. In this context the hitting measure is defined via work of Kaimanovich and Masur \cite{KaimanovichMasur:Poisson}. 
 
Eskin, Mirzakhani and Rafi follow an approach of Connell and Muchnik, who generalized Furstenberg's results to apply, for example, to compact manifolds of negative curvature \cite{ConnellMuchnik:Harmonicity}. Instead of using Brownian motion, this approach pursues a more direct method for solving  Equation \eqref{E:harmonic}.


\bibliography{MM}{}

\newcommand{\etalchar}[1]{$^{#1}$}
\providecommand{\bysame}{\leavevmode\hbox to3em{\hrulefill}\thinspace}
\providecommand{\MR}{\relax\ifhmode\unskip\space\fi MR }
\providecommand{\MRhref}[2]{%
  \href{http://www.ams.org/mathscinet-getitem?mr=#1}{#2}
}
\providecommand{\href}[2]{#2}
\begin{thebibliography}{ABEM12}

\bibitem[ABEM12]{AthreyaBufetovEskinMirzakhani:LatticePoint}
Jayadev Athreya, Alexander Bufetov, Alex Eskin, and Maryam Mirzakhani,
  \emph{Lattice point asymptotics and volume growth on {T}eichm\"{u}ller
  space}, Duke Math. J. \textbf{161} (2012), no.~6, 1055--1111. \MR{2913101}

\bibitem[AGY06]{AvilaGouezelYoccoz:Exponential}
Artur Avila, S\'{e}bastien Gou\"{e}zel, and Jean-Christophe Yoccoz,
  \emph{Exponential mixing for the {T}eichm\"{u}ller flow}, Publ. Math. Inst.
  Hautes \'{E}tudes Sci. (2006), no.~104, 143--211. \MR{2264836}

\bibitem[AH19]{Arana-Herrera:Counting}
Francisco Arana-Herrera, \emph{Counting square-tiled surfaces with prescribed
  real and imaginary foliations and connections to {M}irzakhani's asymptotics
  for simple closed hyperbolic geodesics}, arXiv:1902.05626 (2019).

\bibitem[AN16]{AulicinoNguyen:Three}
David Aulicino and Duc-Manh Nguyen, \emph{Rank two affine manifolds in genus
  3}, arXiv:1612.06970 (2016).

\bibitem[Api18]{Apisa:Hyp}
Paul Apisa, \emph{{${\rm GL}_2\Bbb R$} orbit closures in hyperelliptic
  components of strata}, Duke Math. J. \textbf{167} (2018), no.~4, 679--742.
  \MR{3769676}

\bibitem[BM00]{BekkaMayer:Ergodic}
M.~Bachir Bekka and Matthias Mayer, \emph{Ergodic theory and topological
  dynamics of group actions on homogeneous spaces}, London Mathematical Society
  Lecture Note Series, vol. 269, Cambridge University Press, Cambridge, 2000.
  \MR{1781937}

\bibitem[BM04]{BrooksMakover:Random}
Robert Brooks and Eran Makover, \emph{Random construction of {R}iemann
  surfaces}, J. Differential Geom. \textbf{68} (2004), no.~1, 121--157.
  \MR{2152911}

\bibitem[BQ11]{BenoistQuint:Mesures}
Yves Benoist and Jean-Fran\c{c}ois Quint, \emph{Mesures stationnaires et
  ferm\'{e}s invariants des espaces homog\`enes}, Ann. of Math. (2)
  \textbf{174} (2011), no.~2, 1111--1162. \MR{2831114}

\bibitem[BQ{\etalchar{+}}16]{Translation}
Y~Benoist, J-F Quint, et~al., \emph{Translation of the paper ``{M}esures
  stationnaires et ferm\'{e}s invariants des espaces homog\`enes", by {Y}.
  {B}enoist and {J}.-{F}. {Q}uint}, arXiv:1610.05912 (2016).

\bibitem[BRH17]{BrownRodriguezHertz:Rigidity}
Aaron Brown and Federico Rodriguez~Hertz, \emph{Measure rigidity for random
  dynamics on surfaces and related skew products}, J. Amer. Math. Soc.
  \textbf{30} (2017), no.~4, 1055--1132. \MR{3671937}

\bibitem[BS85]{BirmanSeries:Sparse}
Joan~S. Birman and Caroline Series, \emph{Geodesics with bounded intersection
  number on surfaces are sparsely distributed}, Topology \textbf{24} (1985),
  no.~2, 217--225. \MR{793185}

\bibitem[BS01]{BonahonSozen:Symplectic}
Francis Bonahon and Ya\c{s}ar S\"{o}zen, \emph{The {W}eil-{P}etersson and
  {T}hurston symplectic forms}, Duke Math. J. \textbf{108} (2001), no.~3,
  581--597. \MR{1838662}

\bibitem[BT16]{BridgemanTan:Identity}
Martin Bridgeman and Ser~Peow Tan, \emph{Identities on hyperbolic manifolds},
  Handbook of {T}eichm\"{u}ller theory. {V}ol. {V}, IRMA Lect. Math. Theor.
  Phys., vol.~26, Eur. Math. Soc., Z\"{u}rich, 2016, pp.~19--53. \MR{3497292}

\bibitem[Cal04]{Calta}
Kariane Calta, \emph{Veech surfaces and complete periodicity in genus two}, J.
  Amer. Math. Soc. \textbf{17} (2004), no.~4, 871--908.

\bibitem[CdS01]{daSilva:Lectures}
Ana Cannas~da Silva, \emph{Lectures on symplectic geometry}, Lecture Notes in
  Mathematics, vol. 1764, Springer-Verlag, Berlin, 2001. \MR{1853077}

\bibitem[CM07]{ConnellMuchnik:Harmonicity}
Chris Connell and Roman Muchnik, \emph{Harmonicity of quasiconformal measures
  and {P}oisson boundaries of hyperbolic spaces}, Geom. Funct. Anal.
  \textbf{17} (2007), no.~3, 707--769. \MR{2346273}

\bibitem[CSW]{ChaikaSmillieWeiss}
Jon Chaika, John Smillie, and Barak Weiss, \emph{Tremors and horocycle dynamics
  on the moduli space of translation surfaces}, in preparation.

\bibitem[CW]{ChenWright:WYSIWYG}
Dawei Chen and Alex Wright, \emph{The {WYSIWYG} compactification}, in
  preparation.

\bibitem[DDM14]{DowdallDuchinMasur:Statistical}
Spencer Dowdall, Moon Duchin, and Howard Masur, \emph{Statistical hyperbolicity
  in {T}eichm\"{u}ller space}, Geom. Funct. Anal. \textbf{24} (2014), no.~3,
  748--795. \MR{3213829}

\bibitem[DGZZ]{DelecroixGoujardZografZorich}
Vincent Delecroix, Elise Goujard, Peter Zograf, and Anton Zorich,
  \emph{Square-tiled surfaces of fixed combinatorial type: equidistribution,
  counting, volumes of the ambient strata}, in preparation.

\bibitem[DN09]{DoNorbury:Cone}
Norman Do and Paul Norbury, \emph{Weil-{P}etersson volumes and cone surfaces},
  Geom. Dedicata \textbf{141} (2009), 93--107. \MR{2520065}

\bibitem[Do08]{Do:Thesis}
Norman Nam~Van Do, \emph{Intersection theory on moduli spaces of curves via
  hyperbolic geometry}, Ph.D. thesis, 2008.

\bibitem[Do10]{Do:Asymptoptic}
Norman Do, \emph{The asymptotic weil-petersson form and intersection theory on
  $\mathcal{M}_{g,n}$}, arXiv:1010.4126 (2010).

\bibitem[Do13]{Do:BigSurvey}
\bysame, \emph{Moduli spaces of hyperbolic surfaces and their
  {W}eil-{P}etersson volumes}, Handbook of moduli. {V}ol. {I}, Adv. Lect. Math.
  (ALM), vol.~24, Int. Press, Somerville, MA, 2013, pp.~217--258. \MR{3184165}

\bibitem[Dum15]{Dumas:Skinning}
David Dumas, \emph{Skinning maps are finite-to-one}, Acta Math. \textbf{215}
  (2015), no.~1, 55--126. \MR{3413977}

\bibitem[EFW18]{EskinFilipWright:Hull}
Alex Eskin, Simion Filip, and Alex Wright, \emph{The algebraic hull of the
  {K}ontsevich-{Z}orich cocycle}, Ann. of Math. (2) \textbf{188} (2018), no.~1,
  281--313. \MR{3815463}

\bibitem[Ein06]{Einsiedler:Ratner}
Manfred Einsiedler, \emph{Ratner's theorem on {${\rm SL}(2,\Bbb R)$}-invariant
  measures}, Jahresber. Deutsch. Math.-Verein. \textbf{108} (2006), no.~3,
  143--164. \MR{2265534}

\bibitem[EKL06]{EinsiedlerKatokLindenstrauss:Invariant}
Manfred Einsiedler, Anatole Katok, and Elon Lindenstrauss, \emph{Invariant
  measures and the set of exceptions to {L}ittlewood's conjecture}, Ann. of
  Math. (2) \textbf{164} (2006), no.~2, 513--560. \MR{2247967}

\bibitem[ELa]{EskinLindenstrauss:Random}
Alex Eskin and Elon Lindenstrauss, \emph{{R}andom walks on llocally homogeneous
  spaces}, \url{https://www.math.uchicago.edu/~eskin/RandomWalks/paper.pdf}.

\bibitem[ELb]{EskinLindenstrauss:Zariski}
\bysame, \emph{{Z}ariski dense random walks on homogeneous spaces},
  \url{https://www.math.uchicago.edu/~eskin/RandomWalks/short_paper.pdf}.

\bibitem[ELSV01]{EkedahlLandoShapiroVainshtein:Hurwtiz}
Torsten Ekedahl, Sergei Lando, Michael Shapiro, and Alek Vainshtein,
  \emph{Hurwitz numbers and intersections on moduli spaces of curves}, Invent.
  Math. \textbf{146} (2001), no.~2, 297--327. \MR{1864018}

\bibitem[EM93]{EskinMcMullen:Mixing}
Alex Eskin and Curt McMullen, \emph{Mixing, counting, and equidistribution in
  {L}ie groups}, Duke Math. J. \textbf{71} (1993), no.~1, 181--209.
  \MR{1230290}

\bibitem[EM11]{EskinMirzakhani:Counting}
Alex Eskin and Maryam Mirzakhani, \emph{Counting closed geodesics in moduli
  space}, J. Mod. Dyn. \textbf{5} (2011), no.~1, 71--105. \MR{2787598}

\bibitem[EM18]{EskinMirzakhani:InvariantMeasures}
\bysame, \emph{Invariant and stationary measures for the {${\rm SL}(2,\Bbb R)$}
  action on moduli space}, Publ. Math. Inst. Hautes \'{E}tudes Sci.
  \textbf{127} (2018), 95--324. \MR{3814652}

\bibitem[EMM10]{EskinMirzakhaniMohammadi:Effective}
Alex Eskin, Maryam Mirzakhani, and Amir Mohammadi, \emph{Effective counting of
  simple closed geodesics on hyperbolic surfaces}, arXiv:1905.04435 (2010).

\bibitem[EMM15]{EskinMirzakhaniMohammadi:Isolation}
\bysame, \emph{Isolation, equidistribution, and orbit closures for the {${\rm
  SL}(2,\Bbb R)$} action on moduli space}, Ann. of Math. (2) \textbf{182}
  (2015), no.~2, 673--721. \MR{3418528}

\bibitem[EMMW]{EskinMcMullenMukamelWright:Billiards}
Alex Eskin, Curtis McMullen, Ronen Mukamel, and Alex Wright, \emph{Billiards,
  quadrilaterals and moduli spaces},
  \url{http://www-personal.umich.edu/~alexmw/dm.pdf}.

\bibitem[EMR]{EskinMirzakhaniRafi:Walks}
Alex Eskin, Maryam Mirzakhani, and Kasra Rafi, forthcoming.

\bibitem[EMR19]{EskinMirzakhaniRafi:Strata}
\bysame, \emph{Counting closed geodesics in strata}, Invent. Math. \textbf{215}
  (2019), no.~2, 535--607.

\bibitem[ES16]{ErlandssonSouto:Counting}
Viveka Erlandsson and Juan Souto, \emph{Counting curves in hyperbolic
  surfaces}, Geom. Funct. Anal. \textbf{26} (2016), no.~3, 729--777.
  \MR{3540452}

\bibitem[ES19]{ErlandssonSouto:MCC}
\bysame, \emph{Mirzakhani's {C}urve {C}ounting}, arXiv:1904.05091 (2019).

\bibitem[Esk]{Eskin:Notes}
Alex Eskin, \emph{{L}ectures on the ${SL}(2,\mathbb{R})$ action on moduli
  space}, \url{https://www.math.uchicago.edu/~eskin/luminy2012/lectures.pdf}.

\bibitem[Esk10]{Eskin:Unipotent}
\bysame, \emph{Unipotent flows and applications}, Homogeneous flows, moduli
  spaces and arithmetic, Clay Math. Proc., vol.~10, Amer. Math. Soc.,
  Providence, RI, 2010, pp.~71--129. \MR{2648693}

\bibitem[EU18]{ErlandssonUyanik:Length}
Viveka Erlandsson and Caglar Uyanik, \emph{Length functions on currents and
  applications to dynamics and counting}, arXiv:1803.10801 (2018).

\bibitem[Eyn14]{Eynard:Short}
Bertrand Eynard, \emph{A short overview of the ``{T}opological recursion"},
  arXiv:1412.3286 (2014).

\bibitem[Fil16]{Filip:Splitting}
Simion Filip, \emph{Splitting mixed {H}odge structures over affine invariant
  manifolds}, Ann. of Math. (2) \textbf{183} (2016), no.~2, 681--713.
  \MR{3450485}

\bibitem[FM12]{FarbMargalit:Primer}
Benson Farb and Dan Margalit, \emph{A primer on mapping class groups},
  Princeton Mathematical Series, vol.~49, Princeton University Press,
  Princeton, NJ, 2012. \MR{2850125}

\bibitem[FM14]{ForniMatheus:Intro}
Giovanni Forni and Carlos Matheus, \emph{Introduction to {T}eichm\"{u}ller
  theory and its applications to dynamics of interval exchange transformations,
  flows on surfaces and billiards}, J. Mod. Dyn. \textbf{8} (2014), no.~3-4,
  271--436. \MR{3345837}

\bibitem[For02]{Forni:Deviations}
Giovanni Forni, \emph{Deviation of ergodic averages for area-preserving flows
  on surfaces of higher genus}, Ann. of Math. (2) \textbf{155} (2002), no.~1,
  1--103. \MR{1888794}

\bibitem[FSU18]{FraczekShiUlcigrai:Genericity}
Krzysztof Fraczek, Ronggang Shi, and Corinna Ulcigrai, \emph{Genericity on
  curves and applications: pseudo-integrable billiards, {E}aton lenses and gap
  distributions}, J. Mod. Dyn. \textbf{12} (2018), 55--122. \MR{3808209}

\bibitem[Fu19]{Fu:Cusp}
Ser-Wei Fu, \emph{Cusp excursions for the earthquake flow on the once-punctured
  torus}, arXiv:1506.04534 (2019).

\bibitem[Fur71]{Furstenberg:Walks}
Harry Furstenberg, \emph{Random walks and discrete subgroups of {L}ie groups},
  Advances in {P}robability and {R}elated {T}opics, {V}ol. 1, Dekker, New York,
  1971, pp.~1--63. \MR{0284569}

\bibitem[Gen17]{Gendulphe:Wrong}
Matthieu Gendulphe, \emph{What's wrong with the growth of simple closed
  geodesics on nonorientable hyperbolic surfaces}, arXiv:1706.08798 (2017).

\bibitem[Gha18]{Ghazouani:Tori}
Selim Ghazouani, \emph{Teichm\"uller dynamics, dilation tori and piecewise
  affine circle homeomorphisms}, arXiv:1803.10129 (2018).

\bibitem[Gol84]{Goldman:Symplectic}
William~M. Goldman, \emph{The symplectic nature of fundamental groups of
  surfaces}, Adv. in Math. \textbf{54} (1984), no.~2, 200--225. \MR{762512}

\bibitem[Ham09]{Hamenstadt:Invariant}
Ursula Hamenst\"{a}dt, \emph{Invariant {R}adon measures on measured lamination
  space}, Invent. Math. \textbf{176} (2009), no.~2, 223--273. \MR{2495764}

\bibitem[Ham16]{Hamenstadt:Thin}
Ursula Hamenst{\"a}dt, \emph{Counting periodic orbits in the thin part of
  strata}, preprint (2016).

\bibitem[Hua16]{Huang:Mirzakhani}
Yi~Huang, \emph{Mirzakhani's recursion formula on {W}eil-{P}etersson volume and
  applications}, Handbook of {T}eichm\"{u}ller theory. {V}ol. {VI}, IRMA Lect.
  Math. Theor. Phys., vol.~27, Eur. Math. Soc., Z\"{u}rich, 2016, pp.~95--127.
  \MR{3618187}

\bibitem[Jef05]{Jeffrey:Flat}
Lisa~C. Jeffrey, \emph{Flat connections on oriented 2-manifolds}, Bull. London
  Math. Soc. \textbf{37} (2005), no.~1, 1--14. \MR{2105813}

\bibitem[Ker83]{Kerckhoff:Nielsen}
Steven~P. Kerckhoff, \emph{The {N}ielsen realization problem}, Ann. of Math.
  (2) \textbf{117} (1983), no.~2, 235--265. \MR{690845}

\bibitem[KM96]{KaimanovichMasur:Poisson}
Vadim~A. Kaimanovich and Howard Masur, \emph{The {P}oisson boundary of the
  mapping class group}, Invent. Math. \textbf{125} (1996), no.~2, 221--264.
  \MR{1395719}

\bibitem[Kon92]{Kontsevich:Airy}
Maxim Kontsevich, \emph{Intersection theory on the moduli space of curves and
  the matrix {A}iry function}, Comm. Math. Phys. \textbf{147} (1992), no.~1,
  1--23. \MR{1171758}

\bibitem[LM08]{LindenstraussMirzakhani:PML}
Elon Lindenstrauss and Maryam Mirzakhani, \emph{Ergodic theory of the space of
  measured laminations}, Int. Math. Res. Not. IMRN (2008), no.~4, Art. ID
  rnm126, 49. \MR{2424174}

\bibitem[LMW16]{LelievreMonteilWeiss:Everything}
Samuel Leli\`evre, Thierry Monteil, and Barak Weiss, \emph{Everything is
  illuminated}, Geom. Topol. \textbf{20} (2016), no.~3, 1737--1762.
  \MR{3523067}

\bibitem[LT14]{LuoTan:Dilog}
Feng Luo and Ser~Peow Tan, \emph{A dilogarithm identity on moduli spaces of
  curves}, J. Differential Geom. \textbf{97} (2014), no.~2, 255--274.
  \MR{3263507}

\bibitem[LZ04]{LandoZvonkin:Graphs}
Sergei~K. Lando and Alexander~K. Zvonkin, \emph{Graphs on surfaces and their
  applications}, Encyclopaedia of Mathematical Sciences, vol. 141,
  Springer-Verlag, Berlin, 2004, With an appendix by Don B. Zagier,
  Low-Dimensional Topology, II. \MR{2036721}

\bibitem[Mag17]{Magee:OneSided}
Michael Magee, \emph{Counting one sided simple closed geodesics on {F}uchsian
  thrice punctured projective planes}, arXiv:1705.09377 (2017).

\bibitem[Mar04]{Margulis:Anosov}
Grigoriy~A. Margulis, \emph{On some aspects of the theory of {A}nosov systems},
  Springer Monographs in Mathematics, Springer-Verlag, Berlin, 2004, With a
  survey by Richard Sharp: Periodic orbits of hyperbolic flows, Translated from
  the Russian by Valentina Vladimirovna Szulikowska. \MR{2035655}

\bibitem[Mar17]{Martin:Early}
William~J Martin, \emph{On an early paper of {M}aryam {M}irzakhani},
  arXiv:1709.07540 (2017).

\bibitem[Mas75]{Masur:Geodesics}
Howard Masur, \emph{On a class of geodesics in {T}eichm\"{u}ller space}, Ann.
  of Math. (2) \textbf{102} (1975), no.~2, 205--221. \MR{0385173}

\bibitem[Mas82]{Masur:Interval}
\bysame, \emph{Interval exchange transformations and measured foliations}, Ann.
  of Math. (2) \textbf{115} (1982), no.~1, 169--200. \MR{644018}

\bibitem[Mas85]{Masur:Ergodic}
\bysame, \emph{Ergodic actions of the mapping class group}, Proc. Amer. Math.
  Soc. \textbf{94} (1985), no.~3, 455--459. \MR{787893}

\bibitem[McM07]{McMullen:Dynamics}
Curtis~T. McMullen, \emph{Dynamics of {${\rm SL}_2(\Bbb R)$} over moduli space
  in genus two}, Ann. of Math. (2) \textbf{165} (2007), no.~2, 397--456.
  \MR{2299738}

\bibitem[McM14]{McMullen:ICM}
\bysame, \emph{The work of {M}aryam {M}irzakhani}, Proceedings of the
  {I}nternational {C}ongress of {M}athematicians---{S}eoul 2014. {V}ol. 1,
  Kyung Moon Sa, Seoul, 2014, pp.~73--79. \MR{3728463}

\bibitem[McS98]{McShane:Identity}
Greg McShane, \emph{Simple geodesics and a series constant over {T}eichmuller
  space}, Invent. Math. \textbf{132} (1998), no.~3, 607--632. \MR{1625712}

\bibitem[Mir96]{Mirzakhani:Choosable}
Maryam Mirzakhani, \emph{A small non-{$4$}-choosable planar graph}, Bull. Inst.
  Combin. Appl. \textbf{17} (1996), 15--18. \MR{1386951}

\bibitem[Mir98]{Mirzakhani:Shur}
\bysame, \emph{A simple proof of a theorem of {S}chur}, Amer. Math. Monthly
  \textbf{105} (1998), no.~3, 260--262. \MR{1615548}

\bibitem[Mir04]{Mirzakhani:Thesis}
\bysame, \emph{Simple geodesics on hyperbolic surfaces and the volume of the
  moduli space of curves}, ProQuest LLC, Ann Arbor, MI, 2004, Thesis
  (Ph.D.)--Harvard University. \MR{2705986}

\bibitem[Mir07a]{Mirzakhani:Random}
\bysame, \emph{Random hyperbolic surfaces and measured laminations}, In the
  tradition of {A}hlfors-{B}ers. {IV}, Contemp. Math., vol. 432, Amer. Math.
  Soc., Providence, RI, 2007, pp.~179--198. \MR{2342816}

\bibitem[Mir07b]{Mirzakhani:Invent}
\bysame, \emph{Simple geodesics and {W}eil-{P}etersson volumes of moduli spaces
  of bordered {R}iemann surfaces}, Invent. Math. \textbf{167} (2007), no.~1,
  179--222. \MR{2264808}

\bibitem[Mir07c]{Mirzakhani:Intersection}
\bysame, \emph{Weil-{P}etersson volumes and intersection theory on the moduli
  space of curves}, J. Amer. Math. Soc. \textbf{20} (2007), no.~1, 1--23.
  \MR{2257394}

\bibitem[Mir08a]{Mirzakhani:Earthquake}
\bysame, \emph{Ergodic theory of the earthquake flow}, Int. Math. Res. Not.
  IMRN (2008), no.~3, Art. ID rnm116, 39. \MR{2416997}

\bibitem[Mir08b]{Mirzakhani:Annals}
\bysame, \emph{Growth of the number of simple closed geodesics on hyperbolic
  surfaces}, Ann. of Math. (2) \textbf{168} (2008), no.~1, 97--125.
  \MR{2415399}

\bibitem[Mir10]{Mirzakhani:ICM}
\bysame, \emph{On {W}eil-{P}etersson volumes and geometry of random hyperbolic
  surfaces}, Proceedings of the {I}nternational {C}ongress of {M}athematicians.
  {V}olume {II}, Hindustan Book Agency, New Delhi, 2010, pp.~1126--1145.
  \MR{2827834}

\bibitem[Mir13]{Mirzakhani:Growth}
\bysame, \emph{Growth of {W}eil-{P}etersson volumes and random hyperbolic
  surfaces of large genus}, J. Differential Geom. \textbf{94} (2013), no.~2,
  267--300. \MR{3080483}

\bibitem[Mir16]{Mirzakhani:Orbits}
\bysame, \emph{Counting mapping class group orbits on hyperbolic surfaces},
  arXiv:1601.03342 (2016).

\bibitem[MM95]{MahmoodianMirzakhani:Tripartite}
E.~S. Mahmoodian and Maryam Mirzakhani, \emph{Decomposition of complete
  tripartite graphs into {$5$}-cycles}, Combinatorics advances ({T}ehran,
  1994), Math. Appl., vol. 329, Kluwer Acad. Publ., Dordrecht, 1995,
  pp.~235--241. \MR{1366852}

\bibitem[MMW17]{McMullenMukamelWright:Cubic}
Curtis~T. McMullen, Ronen~E. Mukamel, and Alex Wright, \emph{Cubic curves and
  totally geodesic subvarieties of moduli space}, Ann. of Math. (2)
  \textbf{185} (2017), no.~3, 957--990. \MR{3664815}

\bibitem[Mor05]{Morris:Ratner}
Dave~Witte Morris, \emph{Ratner's theorems on unipotent flows}, Chicago
  Lectures in Mathematics, University of Chicago Press, Chicago, IL, 2005.
  \MR{2158954}

\bibitem[MP17]{MirzakhaniPetri:Lengths}
Maryam Mirzakhani and Bram Petri, \emph{Lengths of closed geodesics on random
  surfaces of large genus}, arXiv:1710.09727 (2017).

\bibitem[MT94]{MargulisTomanov:Invariant}
G.~A. Margulis and G.~M. Tomanov, \emph{Invariant measures for actions of
  unipotent groups over local fields on homogeneous spaces}, Invent. Math.
  \textbf{116} (1994), no.~1-3, 347--392. \MR{1253197}

\bibitem[MV15]{Mirzakhani-Vondrak:SpernerFair}
Maryam Mirzakhani and Jan Vondr\'{a}k, \emph{Sperner's colorings, hypergraph
  labeling problems and fair division}, Proceedings of the {T}wenty-{S}ixth
  {A}nnual {ACM}-{SIAM} {S}ymposium on {D}iscrete {A}lgorithms, SIAM,
  Philadelphia, PA, 2015, pp.~873--886. \MR{3451084}

\bibitem[MV17]{MirzakhaniVondrak:SpernerOptimal}
\bysame, \emph{Sperner's colorings and optimal partitioning of the simplex}, A
  journey through discrete mathematics, Springer, Cham, 2017, pp.~615--631.
  \MR{3726616}

\bibitem[MW02]{Minsky-Weiss:Nondivergence}
Yair Minsky and Barak Weiss, \emph{Nondivergence of horocyclic flows on moduli
  space}, J. Reine Angew. Math. \textbf{552} (2002), 131--177. \MR{1940435}

\bibitem[MW17]{MirzakhaniWright:Boundary}
Maryam Mirzakhani and Alex Wright, \emph{The boundary of an affine invariant
  submanifold}, Invent. Math. \textbf{209} (2017), no.~3, 927--984.
  \MR{3681397}

\bibitem[MW18]{MirzakhaniWright:FullRank}
\bysame, \emph{Full-rank affine invariant submanifolds}, Duke Math. J.
  \textbf{167} (2018), no.~1, 1--40. \MR{3743698}

\bibitem[MZ00]{ManinZograf:Invertible}
Yuri~I. Manin and Peter Zograf, \emph{Invertible cohomological field theories
  and {W}eil-{P}etersson volumes}, Ann. Inst. Fourier (Grenoble) \textbf{50}
  (2000), no.~2, 519--535. \MR{1775360}

\bibitem[MZ15]{MirzakhaniZograf:LargeGenus}
Maryam Mirzakhani and Peter Zograf, \emph{Towards large genus asymptotics of
  intersection numbers on moduli spaces of curves}, Geom. Funct. Anal.
  \textbf{25} (2015), no.~4, 1258--1289. \MR{3385633}

\bibitem[Not18]{NoticesMM}
\emph{Maryam {M}irzakhani: 1977--2017}, Notices Amer. Math. Soc. \textbf{65}
  (2018), no.~10, 1221--1247, H\'{e}l\`ene Barcelo and Stephen Kennedy,
  coordinating editors, With remembrances. \MR{3837071}

\bibitem[PT08]{PapadopoulosTheret:Shift}
Athanase Papadopoulos and Guillaume Th\'{e}ret, \emph{Shift coordinates,
  stretch lines and polyhedral structures for {T}eichm\"{u}ller space},
  Monatsh. Math. \textbf{153} (2008), no.~4, 309--346. \MR{2394553}

\bibitem[Qui16]{Quint:Rigidite}
Jean-Fran\c{c}ois Quint, \emph{Rigidit\'{e} des {${\rm SL}_2(\Bbb R)$}-orbites
  dans les espaces de modules de surfaces plates [d'apr\`es {E}skin,
  {M}irzakhani et {M}ohammadi]}, Ast\'{e}risque (2016), no.~380, S\'{e}minaire
  Bourbaki. Vol. 2014/2015, Exp. No. 1092, 83--138. \MR{3522172}

\bibitem[RS]{RafiSouto:Revisited}
Kasra Rafi and Juan Souto, \emph{Statistics of simple curves on surfaces,
  revisited}, in preparation.

\bibitem[RS17]{RafiSouto:Currents}
\bysame, \emph{Geodesics {C}urrents and {C}ounting {P}roblems},
  arXiv:1709.06834 (2017).

\bibitem[Thu98]{Thurston:Stretch}
William~P Thurston, \emph{Minimal stretch maps between hyperbolic surfaces},
  arXiv preprint math/9801039 (1998).

\bibitem[Vee86]{Veech:Geodesic}
William~A. Veech, \emph{The {T}eichm\"{u}ller geodesic flow}, Ann. of Math. (2)
  \textbf{124} (1986), no.~3, 441--530. \MR{866707}

\bibitem[Ven08]{Venkatesh:EKL}
Akshay Venkatesh, \emph{The work of {E}insiedler, {K}atok and {L}indenstrauss
  on the {L}ittlewood conjecture}, Bull. Amer. Math. Soc. (N.S.) \textbf{45}
  (2008), no.~1, 117--134. \MR{2358379}

\bibitem[Wit91]{Witten:TwoD}
Edward Witten, \emph{Two-dimensional gravity and intersection theory on moduli
  space}, Surveys in differential geometry ({C}ambridge, {MA}, 1990), Lehigh
  Univ., Bethlehem, PA, 1991, pp.~243--310. \MR{1144529}

\bibitem[Wit92]{Witten:Revisited}
\bysame, \emph{Two-dimensional gauge theories revisited}, J. Geom. Phys.
  \textbf{9} (1992), no.~4, 303--368. \MR{1185834}

\bibitem[Wol83]{Wolpert:Homology}
Scott Wolpert, \emph{On the homology of the moduli space of stable curves},
  Ann. of Math. (2) \textbf{118} (1983), no.~3, 491--523. \MR{727702}

\bibitem[Wol07]{Wolpert:Cusps}
Scott~A. Wolpert, \emph{Cusps and the family hyperbolic metric}, Duke Math. J.
  \textbf{138} (2007), no.~3, 423--443. \MR{2322683}

\bibitem[Wol10]{Wolpert:Families}
\bysame, \emph{Families of {R}iemann surfaces and {W}eil-{P}etersson geometry},
  CBMS Regional Conference Series in Mathematics, vol. 113, Published for the
  Conference Board of the Mathematical Sciences, Washington, DC; by the
  American Mathematical Society, Providence, RI, 2010. \MR{2641916}

\bibitem[Wol13]{Wolpert:Mirzakhani}
\bysame, \emph{Mirzakhani's volume recursion and approach for the
  {W}itten-{K}ontsevich theorem on moduli tautological intersection numbers},
  Moduli spaces of {R}iemann surfaces, IAS/Park City Math. Ser., vol.~20, Amer.
  Math. Soc., Providence, RI, 2013, pp.~221--266. \MR{3114687}

\bibitem[Wri15a]{Wright:Cylinder}
Alex Wright, \emph{Cylinder deformations in orbit closures of translation
  surfaces}, Geom. Topol. \textbf{19} (2015), no.~1, 413--438. \MR{3318755}

\bibitem[Wri15b]{Wright:Broad}
\bysame, \emph{Translation surfaces and their orbit closures: an introduction
  for a broad audience}, EMS Surv. Math. Sci. \textbf{2} (2015), no.~1,
  63--108.

\bibitem[Wri16]{Wright:Billiards}
\bysame, \emph{From rational billiards to dynamics on moduli spaces}, Bull.
  Amer. Math. Soc. (N.S.) \textbf{53} (2016), no.~1, 41--56. \MR{3403080}

\bibitem[Wri18]{Wright:Earthquake}
\bysame, \emph{Mirzakhani's work on earthquake flow}, arXiv:1810.07571 (2018).

\bibitem[WX18]{WuXue:Small}
Yunhui Wu and Yuhao Xue, \emph{Small eigenvalues of closed {R}iemann surfaces
  for large genus}, arXiv:1809.07449 (2018).

\bibitem[Zog08]{Zograf:Large}
Peter Zograf, \emph{On the large genus asymptotics of {W}eil-{P}etersson
  volumes}, arXiv:0812.0544 (2008).

\bibitem[Zor14]{Zorich:Magique}
Anton Zorich, \emph{Le th\'{e}or\`eme de la baguette magique de {A}. {E}skin et
  {M}. {M}irzakhani}, Gaz. Math. (2014), no.~142, 39--54. \MR{3278429}

\bibitem[Zor15]{Zorich:Magic}
\bysame, \emph{The magic wand theorem of {A}. {E}skin and {M}. {M}irzakhani},
  arXiv:1502.05654 (2015).

\end{thebibliography}
\footnotesize
\bibliographystyle{amsalpha}

\end{document}